\magnification=1200
\input amstex
\documentstyle{amsppt}
\NoBlackBoxes

 \pageheight{9.1truein}
\def\leftitem#1{\item{\hbox to\parindent{\enspace#1\hfill}}}
\def\Vir{\operatorname{Vir}}
\def\dim{\operatorname{dim}}

\def\a{\alpha}
\def\b{\beta}
\def\f{\forall}
\def\N{\Bbb N}
\TagsOnRight
\def\supp{\operatorname{supp}}
\def\su{\succ}
\def\pr{\prec}

\def\Ind{\operatorname{Ind}}
\def\max{\operatorname{max}}

\def\bC{\Bbb C}
\def\dc{\text{\D c}}
\def\Z{\Bbb Z}
\def\bQ{\Bbb Q}

\def\proclaim#1{\par\vskip.25cm\noindent{\bf#1. }\begingroup\it}
\def\endproclaim{\endgroup\par\vskip.25cm}

\rightheadtext{Classification of weight modules}

\topmatter
\title Classification of irreducible weight
modules over higher rank Virasoro algebras
\endtitle
\author    Rencai Lu, and Kaiming Zhao
\endauthor
\affil Department of mathematics\\ Peiking university\\
       Beijing 100871 \\
       P. R. China\\
       Email: rencail\@yahoo.com.cn\\
       \\and\\
       Department of Mathematics\\ Wilfrid Laurier University\\
Waterloo, Ontario\\ Canada N2L 3C5 \vskip 3pt Institute of Mathematics\\
Academy of Mathematics and System Sciences\\ Chinese Academy of
Sciences\\ Beijing 100080, P.R.
China\\
        Email:  kzhao\@wlu.ca\\
\endaffil
\thanks Research supported by NSERC, and the NSF  of China  (Grants 10371120  and 10431040).
\endthanks

\abstract  Let $G$ be a rank $n$ additive subgroup of $\bC$ and
$\Vir[G]$  the corresponding Virasoro algebra of rank $n$. In the
present paper, irreducible weight modules with finite dimensional
weight spaces over $\Vir[G]$ are completely determined. There are
two different classes of them. One class consists of simple
modules of intermediate series whose weight spaces are all
$1$-dimensional. The other is  constructed by using intermediate
series modules over a Virasoro subalgebra of rank $n-1$. The
classification of such modules over the classical Virasoro algebra
was obtained by O. Mathieu in 1992 using a completely different
approach.

\endabstract

\endtopmatter
\document
\smallskip\bigskip

\smallskip\bigskip
\subhead 1. \ \ Introduction
\endsubhead
\medskip
Let $\bC$ be the  field of complex numbers. The {\bf Virasoro
algebra $\Vir:=\Vir[\Z]$} (over $\bC$) is  the Lie algebra with
the basis $\bigl\{c,d_{i}\bigm|i\in\Z\bigr\}$ and the Lie bracket
defined by
$$\eqalign{ &[c, d_i]=0,\cr &[d_i,
d_j]=(j-i)d_{i+j}+\delta_{i,-j}\frac{i^{3}-i}{12}c,\qquad\forall\,\,i,j\in\Z.\cr}
$$
The structure theory of the Virasoro algebra weight modules with
finite-dimensional weight spaces  is fairly well developed. For
details, we refer the readers to [M], the book [KR] and the
references therein.

\smallskip The centerless Virasoro algebra is actually a Witt
algebra, and generalized Witt algebras in positive characteristic
and characteristic $0$ were studied by many authors, for instance,
Zassenhaus [Z], Kaplansky [K], Ree [R], Wilson [W], Strade [St];
and Osborn [O], Djokovic and Zhao [DZ], Passman [P], Xu [X].

\smallskip Patera and Zassenhaus [PZ] introduced the {\bf generalized
Virasoro algebra} $\Vir[G]$ for any additive subgroup $G$ of
$\bC$. This Lie algebra can be obtained from $\Vir$ by replacing
the index group $\Z$ with $G$ (see Definition 2.1).   If $G\simeq
\Z^n$, then $\Vir[G]$ is called a {\bf rank n Virasoro algebra}
(or a {\bf higher rank Virasoro algebra} if $n\geq 2$).

\smallskip Representations for generalized Virasoro algebras
$\Vir[G]$ have been studied by several authors.  Mazorchuk [Ma1]
proved that all irreducible weight modules with finite dimensional
weight spaces over   $\Vir[\bQ]$ are intermediate series modules
(where $\bQ$ is the field of rational numbers). In [Ma2],
Mazorchuk determined the irreducibility of Verma modules with zero
central charge over higher rank Virasoro algebras. In [HWZ], Hu,
Wang and Zhao obtained a criterion for the irreducibility of Verma
modules over the generalized Virasoro algebra $\Vir[G]$ over an
arbitrary field $F$ of characteristic $0$  ($G$ is an additive
subgroup of $F$). Su and Zhao [SZ] proved that weight modules with
all weight spaces $1$-dimensional are some so-called intermediate
series of modules. In [S1] and [S2], Su proved that the
irreducible weight modules over higher rank Virasoro algebras are
divided into two classes: intermediate series modules, and GHW
modules. In [BZ],  Billig and Zhao constructed a new class of
irreducible weight modules with finite dimensional weight spaces
over some generalized Virasoro algebras.

\vskip 5pt
 The aim of this paper is to complete the classification of
irreducible weight modules with finite dimensional weight spaces
over higher rank Virasoro algebras $\Vir[G]$. The result for $n=1$
was obtained by Mathieu [M] by using a completely different
method.

\vskip 5pt

This paper is arranged as follows.

In Section 2, we collect some known results.
 For any total order $``\su"$ on $G$, which is compatible with
the group addition, and for any $\dc, h\in \bC$, we recall the
definition of the Verma module $ {M}(\dc,h,\su)$ over   $\Vir[G]$
and some known facts about such modules (see [HWZ]). We recall
from [BZ] the construction of a class of irreducible weight
modules with finite dimensional weight spaces over some
generalized (including higher rank) Virasoro algebras. These
modules are denoted by $ V(\a,\b, b, G_0 )$ (see (2.5) for
definition) for some $ \a,\b \in \bC$, $b\in G\setminus\{0\}$, and
a subgroup $G_0$ of $G$ with $G=\Z b \bigoplus G_0$. We also
recall in Theorem 2.5 a useful result from [S1].

In Section 3, we give a classification of  irreducible weight
modules with finite dimensional weight spaces over $\Vir[G]$ for
$G\simeq \Z^n$, i.e., any such   module is either an intermediate
series module $ V'(\a,\b, G)$ or $ V(\a,\b, b, G_0 )$ for suitable
parameters
 (Theorem 3.9). We show that all
GHW modules (see the definition preceding Theorem 2.5) over
$\Vir[G]$ are isomorphic to modules $ V(\a,\b, b, G_0 )$. The main
technique we employ in this paper is to thoroughly study the
weight set  $\supp(V)$ (sometimes also called  the support) of
nontrivial irreducible weight modules $V$ with finite dimensional
weight spaces. We first spend a lot of effort to handle the case
  $n=2$ (Lemma 3.3-Theorem 3.7), and then use induction on $n$ to
deal with all other cases. The induction turns out to be rather
difficult.

 \vskip 5pt
We hope that our results  will  have some applications in physics
since the Lie algebras studied in the present paper have similar
properties as the classical Virasoro algebra which is widely used
in physics.

 \vskip 5pt {\bf Acknowledgement.} The authors would like
to thank sincerely Prof.\ Volodymyr Mazorchuk and the referee for
scrutinizing all steps in the old version of the paper, pointing
out several inaccuracies, and making valuable suggestions to
improve the paper. The second author likes to express his thanks
to Prof. D. Z. Djokovic for helping with the final version.

\smallskip\bigskip
\subhead 2. \ \ Weight modules over Generalized Virasoro algebras
\endsubhead
\medskip

In this section we recall the construction of various modules and
collect some known results for later use.

\proclaim{Definition 2.1} Let $G$ be a nonzero additive subgroup
of $\bC$. The {\bf generalized Virasoro algebra $\Vir[G]$} (over
$\bC$) is the Lie algebra with the
 basis $\bigl\{c,d_{x}\bigm|x\in G\bigr\}$ and
the Lie bracket defined by $$\eqalign{ &[c, d_{x}]=0,\cr &[d_{x},
d_y]=(x-y)d_{x+y}+\delta_{x,-y}\frac{x^{3}-x}{12}c,\qquad
\forall\,\,x,y\in G.\cr} $$
\endproclaim

It is clear that $\Vir[G]\simeq \Vir[aG]$ for any $a\in \bC^*$.
 For any $x\in G^{\ast}:= G\setminus\{0\}$,
$\Vir[x\Z]$ is a Lie subalgebra of $\Vir[G]$ isomorphic to $\Vir$.

Fix a total order ``$\su$'' on $G$ which is compatible with the
addition, i.e., $x\su y$ implies $x+z\su y+z$ for any $z\in G$.
Let
$$ G_{+}:=\bigl\{x\in G\bigm|x\su 0\bigr\},\quad
G_{-}:=\bigl\{x\in G\bigm|x\pr 0\bigr\}. $$ Then
$G=G_{+}\cup\bigl\{0\bigr\}\cup G_{-}$ and we have the triangular
decomposition $$\Vir[G]=\Vir[G]_+\oplus\Vir[G]_-\oplus\Vir[G]_0,$$
where $\Vir[G]_+=\oplus_{x\in G_+}\bC d_x$,
$\Vir[G]_-=\oplus_{x\in G_-}\bC d_x$, $\Vir[G]_0=\bC d_0+\bC c$.

It is clear  that either
$$\#\bigl\{y\in G\bigm|0\pr y\pr
x\bigr\}=\infty\,\,\,\, \forall\,\,x\in G_{+},\,\, \tag2.1$$ or
$$ \exists \,\,a\in G_{+},  \,\, \#\bigl\{y\in
G\bigm|0\pr y\pr a\bigr\} =0.
 \tag2.2$$ We say that the order is
{\it dense} resp. {\it discrete} if (2.1) resp. (2.2) holds.

\vskip 5pt A ${\Vir}[G]$-module $V$ is called {\bf trivial} if
${\Vir}[G]V=0$. For any ${\Vir}[G]$-module $V$ and $\dc,\lambda\in
\bC$, let $V_{\dc,\lambda}:=\bigl\{v\in V\bigm|d_{0}v=\lambda v,
cv=\dc v\bigr\}$ denote the {\bf weight space} of $V$
corresponding to a weight $(\dc,\lambda)$. When $c$ acts as the
scalar $\dc$ on the whole space $V$, we shall simply write
$V_{\lambda}$ instead of $V_{\dc,\lambda}$.

A $\Vir[G]$-module $V$ is called a {\bf weight module} if $V$ is
the sum of its weight spaces. For a weight module $V$ we define
$\supp(V):=\bigl\{\lambda\in \bC\bigm|V_{\lambda}\neq 0\bigr\}$,
which is generally called the {\bf weight set} (or the {\bf
support}) of $V$.

\vskip 5pt For any Lie algebra $L$, we shall use $U(L)$ to denote
its universal enveloping algebra. For any $\dc, h\in \bC$, let
$I(\dc,h,\su)$ be the left ideal of $U:=U(\Vir[G])$ generated by
the elements
$$ \bigl\{d_{i}\bigm|i\in G_+ \bigr\}\bigcup\bigl\{d_0-h\cdot 1,
c-\dc\cdot 1\bigr\}. $$ Then the {\bf Verma module} with the
highest weight $(\dc, h)$ for $\Vir[G]$ is defined as
$M(\dc,h,\su):=U/I(\dc,h,\su)$. This module has a basis consisting
of the following vectors
$$ d_{-i_1}d_{-i_2}\cdots d_{-i_k}v_{h},\quad
\text{$k\in{\N}\cup\{0\}$, $i_{j}\in G_+,\, \forall\,\,j$ and
$i_{k}\geq\cdots\geq i_2\geq i_1>0,$} $$ where
$v_h=1+I(\dc,h,\su)$ is the highest weight vector. Let
$V(\dc,h,\su)$ be the unique irreducible quotient of
$M(\dc,h,\su)$. Let us recall

\proclaim{Theorem 2.2 ([HWZ, Theorem 3.1])} Let $\dc, h\in \bC $.

1) Assume that the order ``$\su$"   is dense. Then the Verma
module $M(\dc,h,\su)$ is an irreducible $\Vir[G]$-module if and
only if $(\dc,h)\neq (0,0)$. Moreover,  $$
M'(0,0,\su):=\sum_{i_{1},\cdots,i_{k}\in G_{+},k>0} \bC
d_{-i_1}\cdots d_{-i_k}v_0$$ is an irreducible submodule of
$M(0,0,\su)$.

2) Assume that the order ``$\su$"  is discrete. Then the Verma
module $M(\dc,h,\su)$ is an irreducible $\Vir[G]$-module if and
only if for the minimal positive element $a\in G$ with respect to
``$\su$", the $\Vir[a\Z]$-module $M_{a}(\dc,h,\su)=U(\Vir[\Z
a])v_h$ is irreducible.
\endproclaim

Now we give another class of   $\Vir[G]$-modules $V(\a,\b,G)$ for
any $\a,\b\in \bC$ (see [SZ]). These $\Vir[G]$-modules all have
basis $\{v_x\,|\,x\in G\}$ with actions given by the following
formula
$$ c v_y=0,\,\,\,\, d_x v_y=(\a+y+x \b)v_{x+y},\,\,\forall\,\,x,y\in
G.$$

One knows from [SZ] that $ V(\a,\b,G)$ is reducible if and only if
$\a\in G$ and $\b\in\{0,1\}$. By $V'(\a,\b,G)$ we denote the
unique nontrivial irreducible sub-quotient of $ V(\a,\b,G)$. Then
$\supp(V'(\a,\b,G))=\a+G$ or $\supp(V'(\a,\b,G))=G\setminus\{0\}$.
We now recall

\vskip 5pt \noindent {\bf  Theorem 2.3 ([SZ, Theorem 4.6]).} {\it
Let $V$ be a nontrivial irreducible weight module over $\Vir[G]$
with all weight spaces $1$-dimensional. Then $V \simeq V'(a,b,G)$
for some $a,b\in\bC$.}

\vskip 5pt
 Now we assume that $G= \Z b\oplus G_0\subset\bC$
where $0\neq b \in \bC$ and $G_0$ is a nonzero subgroup of $\bC$.
(Note that some $G$ lack this property). We temporarily set
$L=\Vir[G]$. For any $i\in\Z$, we set
$$L_{i b}=\oplus_{a\in G_0}\bC d_{i b+a},$$
$$L_+=\oplus_{i>0} L_{i b},\,\, L_-=\oplus_{i<0} L_{i b},\,\,\,
L_0\simeq \Vir[G_0].$$

For any $\a,\b\in\bC$, we have the irreducible $L_0$-module
$V'(\a,\b,G_0)$. We extend the $L_0$-module structure on
$V'(\a,\b,G_0)$ to an $(L_++L_0)$-module structure by defining
$L_+V'(\a,\b,G_0)=0$. Then we obtain the induced $L$-module
$$\aligned \bar M(b,G_0,V'(\a,\b,G_0))&=\Ind_{L_++L_0}^{L}V'(\a,\b,G_0)\\
=&U(L)\otimes_{U(L_++L_0)}V'(\a,\b,G_0).\endaligned \tag2.3$$

 As vector spaces,
 $ \bar M(b,G_0,V'(\a,\b,G_0))\simeq U(L_-)\otimes_{\bC}V'(\a,\b,G_0)$.
  The $L$-module   $\bar M(b,G_0,V'(\a,\b,G_0))$
 has a unique maximal proper submodule $J$. Then we obtain the irreducible quotient
 module
 $$ {M}( b,G_0,V'(\a,\b,G_0))= \bar M(b,G_0,V'(\a,\b,G_0))/J.\tag 2.4 $$
 It is clear that this module is uniquely determined by
 $\a,\b, b$ and  $G_0$; and that $$\supp({M}( b,G_0,V'(\a,\b,G_0)))=
 \Z^+b+G_0\,\,\,\roman{or}\,\,\,(\Z^+b+G_0)\setminus\{0\}.\tag2.5$$
 Note that $b$ can be replaced by any element in $b+G_0$.

To simplify notation, set
$$V= V(\a,\b, b,G_0)={M}(b,G_0,V'(\a,\b,G_0)). \tag 2.6$$
It is clear that $ V= \oplus_{i\in\Z_+} V_{-i b+\a+G_0}$, where
 $$ V_{-i b+\a+G_0}=\oplus_{a\in G_0} V_{-i b+\a+a},\,\,\,
  V_{-i b+\a+a}=\{v\in   V\,\,|\,\,d_0v=(-i b+\a+a)v\}.$$

Now we recall

 \medskip
\noindent {\bf   Theorem 2.4 ([BZ, Theorem 3.1]).}   {\it All
weight spaces of the $\Vir[G]$-module $V(\a,\b, b,G_0)$, defined
above, are finite dimensional. More precisely, $\dim  V_{-i
b+\a+a}\break \le (2i+1)!!$ for all $ i\in \N,\, a\in G_0.$}

\medskip

{\bf From now on in this paper} we assume that $G\simeq \Z^n$ for
some integer $n>1$, $V=\oplus_{x\in G} V_{a+x}$ is an irreducible
weight module over $\Vir[G]$ with finite dimensional weight spaces
(i.e., $\dim V_{a+x}<\infty$ for all $x\in\bC$) where $a\in\bC$.
If there exists $N\in \N$ such that $\dim V_{a+x}<N$ for all $x\in
\bC$, we say that $V$ is {\bf uniformly bounded}. If there exists
a $\Z$-basis $B=\{b_1,\cdots  ,b_n\}$ of $G$ and $v_{\Lambda_0}\in
V_{\Lambda_0}$ such that
$$d_xv_{\Lambda_0}=0, \f \,\,0\ne x\in \Z^+b_1+\cdots  +\Z^+b_n,$$
we say that $V$ is a {\bf generalized highest weight module} (GHW
module for short) with GHW $\Lambda_0$ w.r.t. $B$ (see [S1]). The
vector $v_{\Lambda_0}$ is called a GHW vector with respect to $B$,
or simply a GHW vector. Finally we recall

 \medskip
\noindent {\bf  Theorem 2.5 ([S1, Theorem 1.2]).}   {\it Suppose
that $G\simeq \Z^n$, $n>1$ and $V$ is a nontrivial irreducible
weight $\Vir[G]$-module with finite dimensional weight spaces.

(a) If $V$ is uniformly bounded,  then $V\simeq V'(\a,\b,G)$ for
suitable $\a,\b\in\bC$.

(b) If $V$ is not uniformly bounded, then $V$ is a GHW module.}

\smallskip\bigskip
\subhead 3. \ \ Classification of weight modules
\endsubhead
\medskip

In this section we give a classification of all irreducible weight
modules with finite dimensional weight spaces over higher rank
Virasoro algebras. More precisely, we prove that any such module
is either $ V'(\a,\b, G)$ or $ V(\a,\b, b, G_0 )$ (Theorem 3.9).
To this end, by Theorem 2.5, we need only  study GHW modules.

Recall that $G$ is an additive subgroup of $\bC$ with $G\simeq
\Z^n$ and $n>1$, and that $V=\oplus_{x\in G} V_{a+x}$ is an
irreducible weight module over $\Vir[G]$ with finite dimensional
weight spaces.

By ``$\su"$ we denote the  lexicographic order on $\Z^n$, i.e.,
$(x_1,\cdots  ,x_n)\su (y_1,\cdots  ,y_n)$ if and only if there
exists $s$: $1\leq s \leq n $ such that $x_i=y_i$ for $1\leq i
\leq s-1$ and $x_s>y_s$.

We write $(x_1,\cdots  ,x_n)> (y_1,\cdots  ,y_n)$ if $x_i>y_i$ for
$1\leq i \leq n$; and $(x_1,\cdots  ,x_n)\ge (y_1,\cdots ,y_n)$ if
$x_i\ge y_i$ for $1\leq i \leq n$.

  In this section, the letters
$i,j,k,l,m,n,p,q,r,s,t,x,y$   denote integers. For convenience, we
set $[p,q]=\{x| x \in \Z, p \leq x \leq q\}$ and define similarly
the infinite intervals $(-\infty, p]$, $[q,\infty)$ and $(-\infty,
+\infty)$. For  $a \in G$ or $ S\subset G$, we denote by $\Vir[a]$
or $\Vir[S]$ the subalgebra of $\Vir[G]$ generated by $\{d_{\pm
a}, d_{\pm 2a}\}$ or $\{d_{\pm a}, d_{\pm 2a}| a\in S\}$,
respectively.

\proclaim{Lemma 3.1}  Suppose that $B=({b_1},{b_2},\cdots  ,b_n)$
is a $\Z$-basis of $G$  and $n\geq 2$. Let $V$ be a nontrivial
irreducible GHW $\Vir[G]$-module with GHW $\Lambda_0$ w.r.t. $B$.

 a)  For any $v \in V$, there exists $p> 0$ such that
 $d_{i_1{b_1}+i_2{b_2}+\cdots  +i_nb_n} v = 0$ for all $(i_1,i_2,\cdots  ,i_n)
  \geq (p,p,\cdots  ,p)$.

 b) If $\Lambda_0+i_1b_1+i_2b_2+\cdots  +i_nb_n \in \supp (V)$, then for any
  positive integers $k_1,k_2,\cdots  ,k_n$, there exists
  $m\geq 0$ such that  $\{x \in \Z\,\,\, |\,\,\,\Lambda_0+
  i_1b_1+i_2b_2+\cdots  +i_nb_n+x(k_1b_1+k_2b_2+\cdots  +k_nb_n) \in
 \supp (V)\}=(-\infty, m]$.

 c) Let S be any subgroup of G of rank n, then any nonzero
 $\Vir[S]$-submodule of $V$ is nontrivial.

 d) There exists a $\Z$-basis $B'=\{b_1',b_2',\cdots  ,b_n'\}$ of G such
 that

 \hskip 10pt d1) V is a GHW module with GHW   $\Lambda_0$ w.r.t. $B'$;

\hskip 10pt  d2) $(\Lambda_0+\Z^+ b_1'+\Z^+ b_2'+\cdots  +
\Z^+b_n')\bigcap
 \supp (V)=\{\Lambda_0\}$;

\hskip 10pt  d3) $(\Lambda_0-\Z^+ b_1'-\Z^+ b_2'-\cdots  -
\Z^+b_n') \bigcap
 \supp (V)=\Lambda_0-\Z^+ b_1'-\Z^+ b_2'-\cdots  -
\Z^+b_n'$;

\hskip 10pt  d4) $\Lambda_0+k_1b_1'+k_2 b_2'+\cdots  +
k_nb_n'\notin \supp (V),\,\,\forall\,\,(k_1,k_2,\cdots  ,k_n) \geq
 (i_1,i_2,\cdots  ,i_n)$

 \hskip 1.5cm if  $\Lambda_0+i_1b_1'+i_2b_2'+\cdots  +i_n b_n'\notin \supp (V)$;

\hskip 10pt  d5) $\Lambda_0+k_1b_1'+k_2 b_2'+\cdots  + k_nb_n' \in
\supp (V),\,\,\forall\,\,(k_1,k_2,\cdots  ,k_n)\leq
 (i_1,i_2,\cdots  ,i_n)$

 \hskip 1.5cm if $\Lambda_0+i_1b_1'+i_2b_2'+\cdots  +i_n b_n'\in \supp (V)$;

\hskip 10pt  d6) For any
 $0\ne (k_1,k_2,\cdots  ,k_n) \ge 0$ and $(i_1,i_2,\cdots  ,i_n)\in\Z^n$, we have

  \hskip 1.5cm  $\{ x \in \Z |
 \Lambda_0+\sum_{l=1}^ni_l b_l' + x( \sum_{l=1}^nk_lb_l')
 \in \supp (V)\}= (-\infty, m]$

 \hskip 1.5cm  for some $m \in \Z$.
\endproclaim
\demo{Proof}  For $n=2$ a slightly weaker form of this lemma  is a
combination of several lemmas in [S2].

 a) Without loss of generality, we may assume that $v=u
 v_{\Lambda_0}$, where
 $$u=d_{i_1^{(1)}b_1+i_2^{(1)}b_2+\cdots  +i_n^{(1)}b_n}
 d_{i_1^{(2)}b_1+i_2^{(2)}b_2+\cdots  +i_n^{(2)}b_n}
\cdots  d_{i_1^{(m)}b_1+i_2^{(m)}b_2+\cdots  +i_n^{(m)}b_n}$$ $$
\hskip -8cm \in U(\Vir[G]).$$ Take $p=\max\{-\sum_{i_{1}^{(s)}<0}
i_{1}^{(s)}, -\sum_{i_{2}^{(s)}<0} i_{2}^{(s)},\cdots
,-\sum_{i_{m}^{(s)}<0} i_{m}^{(s)}\}+1$. By induction on $m$, and
using the Lie bracket in $\Vir[G]$, we easily obtain
$$d_{i_1{b_1}+i_2{b_2}+\cdots  +i_nb_n} v =
0,\,\,\,\forall\,\,\,(i_1,i_2,\cdots  ,i_n) \geq (p,p,\cdots
,p).$$

\medskip
 b) Let
 $J= \{x \in \Z\,\,\, |\,\,\,\Lambda_0+\sum_{l=1}^ni_lb_l+x(\sum_{l=1}^nk_lb_l
 ) \in
 \supp (V)\}$.
 \medskip
 {\bf Claim 1.}  {\it For any nonzero $ v \in V$, we have
 $d_{-(k_1b_1+k_2b_2+\cdots  +k_nb_n)}v \neq
 0.$}
\medskip
  {\it Proof of Claim 1.} Suppose that
  $d_{-(k_1b_1+k_2b_2+\cdots  +k_nb_n)}v  =0$ for some nonzero $ v \in V$.
    Let $p$ be as in a). Then
  $d_{-(k_1b_1+k_2b_2+\cdots  +k_nb_n)}$ and
  $d_{b_i+p(k_1b_1+k_2b_2+\cdots  +k_nb_n)}$ for $i\in[1,n]$ act trivially
  on $v$. Since  $\Vir[G]$ is generated by these elements, we see that $\Vir[G]v=0$,
   contradicting  the fact that $V$ is a nontrivial irreducible
module. Claim 1 follows.
  \medskip

  It follows from this claim that $J=(-\infty ,m]$ for some $m\geq 0$ or $J=\Z$.

   Suppose that $J=\Z$. For any $x\in\Z$,
   let
   $$\lambda_x=\Lambda_0+i_1b_1+i_2b_2+\cdots  +i_nb_n+x(k_1b_1+k_2b_2+\cdots  +k_nb_n).$$
We know that
    $\Vir^{[k]}:=\Vir[k_1b_1+k_2b_2+\cdots  +k_nb_n]$ is a rank one
   Virasoro subalgebra, and $W=\oplus_{x\in\Z}
   V_{\lambda_x}$
   is a $\Vir^{[k]}$-module.
    From a) and a well-known result in [M, Lemma 1.6] for any $x \in \Z$ there exists $y
   \geq x $ such that
   $V_{\lambda_{y}}$
contains a $\Vir^{[k]}$ primitive vector (a nonzero weight vector
$v$ such that $d_{l(k_1b_1+k_2b_2+\cdots  +k_nb_n)}v=0$ for all
$l\in\N$). So there are infinitely many nontrivial highest weight
$\Vir^{[k]}$-modules  having the same weight $\lambda_0$, which
implies $\dim V_{\lambda_0}=\infty$. This contradiction yields
that $J\ne\Z$. Hence   b) is proved.
\medskip

 c)
 For any $p$, let $I_1=(p+1,p,\cdots  , p)$,
 $I_2=(p+2,p+1,p,p,\cdots  ,p)$,
 $I_k=I_1+(0,0,\delta_{3,k},\cdots  , \delta_{n,k}) \in \Z^n, k=3,\cdots  ,n$.
 Let $$A=\left(\matrix I_1\crcr
I_2\crcr I_3\crcr \vdots\crcr I_n \endmatrix\right). \tag 3.1$$
Then $\det(A)=1$. Suppose that there exists a rank $n$ subgroup
$S$ of $G$ and a nonzero $v_0 \in V$ such that $\Vir[S]v_0=0$. Now
take $p$ as in a), that is, $d_{i_1{b_1}+i_2{b_2}+\cdots  +i_nb_n}
v_0 = 0$ for all $(i_1,i_2,\cdots  ,i_n) \geq (p,p,\cdots  ,p)$.
Let $(b_1^*,b_2^*,\cdots  ,b_n^*)=(b_1,b_2,\cdots  ,b_n)A$. Then
$d_{b_i^*}v_0=0$ for all $i=1,2,\cdots  ,n$. Since $G/S$ is a
finite group, there exists some $i>0$, such that
$-i(b_1^*+b_2^*+\cdots  +b_n^*) \in S$. Clearly $d_{-b_1^*},
d_{-b_2^*},\cdots  ,d_{-b_n^*}$ belongs to the subalgebra
generated by the elements: $d_{-i(b_1^*+b_2^*+\cdots  +b_n^*)},
d_{b_i^*},i=1,2,\cdots  ,n$; and   $\Vir[G]$ is generated by
$d_{\pm b_i^*},i=1,2,..,n$. Hence we have $\Vir[G]v_0=0$, a
contradiction to the fact that $V$ is nontrivial. \medskip

d) By b) we can suppose that $\{x \in \Z\,\,\,
|\,\,\,\Lambda_0+x(b_1+b_2+\cdots  +b_n) \in
 \supp (V)\}=(\infty, p-2]$ for some $p\geq 2$. Take $A$ as in (3.1), and
 $(b_1',b_2',\cdots  ,b_n')=(b_1,b_2,\cdots  ,b_n)A$. One can easily check
 d1)-d6) by using b) and Claim 1. We omit
 the details.
\hfill $ $\qed \enddemo

To better understand the proof of Lemma 3.1 d) and the lemmas that
follow it might help if one draws a diagram in the $Ob_1b_2$-plane
for $n=2$ to describe those sets. For instance, if
$\lambda=x_1b_1+x_2b_2$ in the first quadrant, i.e.,
$x_1>0,x_2>0$, then $\Lambda_0+\lambda\notin\supp(V)$ and
$\Lambda_0-\lambda\in\supp(V)$.

In the next lemma we do not assume the irreducibility of $V$.

 \proclaim{Lemma 3.2}
If $V$ is a nonzero uniformly bounded weight module over
$\Vir[G]$, then $V$ has an irreducible submodule.
\endproclaim
\demo{Proof}
 Fix $a\in\supp(V)$. Then $\oplus_{g\in G }V_{g+a}$
is a $\Vir[G]$-submodule.  Thus it is enough to prove the lemma
for $V=\oplus_{g\in G }V_{g+a}$. We may  assume that $V$ does not
have any nonzero trivial submodules. So we can further assume that
$a\ne0$.

We shall prove the lemma by induction on $\dim V_a$.

If $\dim V_a=1$, let $W$ be the submodule generated by $V_a$. We
know that there exists a maximal proper submodule $W'$ of $W$ not
containing $V_a$. So $W/W'$ is irreducible. By Theorem 2.5, we
know that $W/W'$ is a $V'(\a,\b,G)$. So $W'_a=0$. If $W'_{a'}\ne0$
for $a'\ne a$, then consider the $\Vir[a-a']$-module generated by
$W'_{a'}$. From the Virasoro algebra theory, we see that $a'=0$,
hence $\Vir[G] W'=0$, a contradiction. So $W'=0$ and $W$ itself is
an irreducible $\Vir[G]$-submodule.
  The lemma follows in this case.

In general, for any nonzero $v\in V_a$, let $W$ be the submodule
generated by $v$. By Zorn's Lemma, there exists a maximal proper
submodule $W'$ of $W$  not containing $v$. By Theorem 2.5,
$W/W'\simeq V'(\a,\b,G)$ for some $\a,\b\in\bC$. If $W'=0$ we are
done. If $W'\ne0$, then, applying the inductive hypothesis to
$W'$, we have  an irreducible submodule of $W'$. The lemma is
proved. \hfill $ $\qed
\enddemo

In the rest of this section we further assume that $V=\oplus_{g\in
G} V_{\Lambda_0+g}$ is a nontrivial irreducible GHW
$\Vir[G]$-module with GHW $\Lambda_0$ w.r.t. $B=\{b_1,b_2,\cdots
,b_n\}$, where $\Lambda_0\in \bC$, and  $B$ satisfies the
properties of Lemma 3.1 d).

\proclaim{Lemma 3.3}  If there exist $(i_1,i_2,\cdots  ,i_n),
(k_1,k_2,\cdots  ,k_n) \in \Z^n$ with $k_1,\cdots  ,k_n $
relatively prime, and $(s_1,\cdots  ,s_n)>0$ satisfying $$\{
\Lambda_0+\sum_{t=1}^ni_t b_t+ \sum_{t=1}^n x_ts_tb_t |
(x_1,x_2,\cdots  ,x_n)\break  \in \Z^n , \sum_{t=1}^n
k_ts_tx_t=0\}
 \bigcap \supp (V)=\emptyset,$$ then
  $V\simeq   M(b', G_0, V'(\a,\b, G_0))$
 for some $ \a,\b \in \bC$, and $G=\Z b' \bigoplus G_0$,
 where $0 \neq b' \in \bC, G_0$
 is a subgroup of $G$.
\endproclaim
{\bf Remark.} The above condition  means that a lattice in some
affine hyperplane of $\Z^n$ orthogonal to $(k_1,k_2,\cdots,k_n)$
contains no weights of $V$.

\demo{Proof} As mentioned earlier, to understand the proof of this
lemma better it may be helpful to sketch in the $Ob_1b_2$-plane
for $n=2$ the sets used in the proof.

By Lemma 3.1 d6), we have $k_i> 0$ for all $i=1,2,\cdots  ,n$ or
$k_i< 0$ for all
 $i=1,2,\cdots  ,n$. We may assume that $(k_1,k_2,\cdots  ,k_n)> 0$.
 Let $$G_0=\{\sum_{t=1}^nx_t b_t
 \in G\,\,|\,\, \sum_{i=1}^n k_ix_i=0
 \}.\tag3.2
$$

 {\bf Claim 1.}  {\it There exists $m_0 \in
    \Z$ such that}
    $$\{  \Lambda_0+\sum_{t=1}^nx_t b_t | (x_1,x_2,\cdots  ,x_n) \in \Z^n ,
 \sum_{i=1}^n k_ix_i \geq m_0\}\bigcap \supp(V)=\emptyset.\tag3.3$$

 {\it Proof of Claim 1.}  Let
$A_t=s_ts_{1}(-\delta_{1,t}k_{1}+k_t,-\delta_{2,t}k_{1} ,\cdots
,-\delta_{n,t}k_{1})$ whose corresponding element in $G$ is
$s_ts_{1}(-k_{1}b_t+k_tb_{1})\in G_0$. Note that $k_{1}\neq 0$.
One may easily check that for any $(z_1,z_2,\cdots  ,z_n) \in
\Z^n$ with $\sum_{t=1}^n z_tk_t \geq 0$, there exist suitable $l_t
\in \Z,t=1,2,\cdots  ,n$, such that
$$
(z_1,z_2,\cdots  ,z_n)=(z_1',z_2',\cdots  ,z_n')+\sum_{t=1}^n
l_tA_t, \tag 3.4
$$
where $0\geq z_t'> -k_{1}s_{1}s_t$ for all $t \in \{2,3,\cdots
,n\}$. Hence $z_{1}'\geq 0$. Now let
$N=\max\{k_{1}s_{1}s_1,k_{1}s_{1}s_2,\cdots  ,k_{1}s_{1}s_n\}$,
and $m_0=\sum_{t=1}^n k_t(N+i_t)$. Then using (3.4), for any
$(x_1,x_2,\cdots  ,x_n) \in \Z^n$ with $ \sum_{i=1}^n k_ix_i \geq
m_0$ we have
$$(x_1,x_2,\cdots  ,x_n)-(i_1+N,i_2+N,\cdots  ,i_n+N)=(x_1',x_2',\cdots  ,x_n')+\sum_{t=1}^n
l_tA_t,$$ where $0\geq x_t'> -k_{1}s_{1}s_t$, i.e.,
$$(x_1,x_2,\cdots  ,x_n)=(i_1,i_2,\cdots  ,i_n)+(N+x_1',N+x_2',\cdots  ,N+x_n')+\sum_{t=1}^n
l_tA_t.\tag3.5$$ Let $(y_1,y_2,\cdots  ,y_n)= \sum_{t=1}^n
l_tA_t$. We have
$$\Lambda_0+\sum_{t=1}^nx_t b_t=
\Lambda_0+\sum_{t=1}^ni_t b_t+
\sum_{t=1}^ny_tb_t+\sum_{t=1}^n\Z^+b_t. $$ Note that
$\sum_{t=1}^ny_tb_t=\sum_{t=1}^ny'_ts_tb_t$ with
$\sum_{t=1}^ny'_ts_tk_t=0$. From the assumption we know that
$$\Lambda_0+\sum_{t=1}^ni_t b_t+
\sum_{t=1}^ny_tb_t \notin \supp(V).\tag3.6$$ By applying Lemma 3.1
d4) we obtain $\Lambda_0+\sum_{t=1}^nx_t b_t\notin  \supp(V)$. The
claim follows.\qed

\vskip 5pt
 From Claim 1 we have a unique integer $m$
 with the following two properties:

 1) $\{ \Lambda_0+\sum_{t=1}^nx_t b_t \in \supp(V)\,\,
 |\,\, x_1,x_2,\cdots  ,x_n \in \Z, \sum_{i=1}^n k_ix_i \geq m\}
 =\emptyset$, and

 2) $P:=\{  \Lambda_0+\sum_{t=1}^nx_t b_t
  \in \supp (V)\,\,| \,\,x_t \in \Z, \sum_{i=1}^n k_ix_i =m-1\}\neq
 \emptyset$.

\vskip 5pt
 Fix some $b_1'=t_1b_1+t_2b_2+\cdots  +t_nb_n$ with
 $\sum_{i=1}^n k_it_i=1$.  Since for any $g=\sum_{i=1}^ng_ib_i\in G$,
 we see that $g-(\sum_{i=1}^n k_ig_i)b'_1\in G_0$, then  $G=\Z b_1'\oplus G_0$.
 Fix  $\lambda_0 \in  P  $.
We have $ P =(\lambda_0+G_0)\bigcap \supp (V)$.
 Let W=$\bigoplus_{ \lambda\in \lambda_0+G_0} V_{\lambda}$,
 which is a $\Vir[G_0]$-submodule of $V$.

\vskip 5pt
 {\bf Claim 2.} {\it W is a uniformly bounded  $\Vir[G_0]$-module.}
\vskip 5pt {\it Proof of Claim 2.}
 Let $0\neq w \in V_{\lambda}$ for some $\lambda \in  P $.
 Noting that $ (P +G_0+b_1')\bigcap \supp (V)=\emptyset$,
 and that for any $a_0 \in G_0$, the set $\{d_{a+b_1'}, d_{-a_0-b_1'}|a \in
 G_0\}$ generates the Lie algebra $\Vir[b_1',G_0]=\Vir[G]$, we deduce $d_{-a_0-b_1'} w
 \neq 0$ for any $a_0 \in G_0$.  Thus we obtain a linear
  injection
$d_{-a_0-b_1'}:\,\,V_{\lambda+a_0} \longrightarrow V_ {\lambda-
b_1'}$. Thus $\dim V_{\lambda+a_0} \leq \dim V_{\lambda-b_1'}$ for
all $a_0 \in G_0$, i.e., $W$ is uniformly bounded. Claim 2
follows. \qed

\vskip 5pt
 By Lemma 3.2, $W$ has  an irreducible
 $\Vir[G_0]$-submodule $W'$. By Theorem 2.5, any irreducible uniformly bounded module
 is either trivial or isomorphic to $V'(\a,\b,G_0)$ for some $(\a,\b)\in
 \bC^2$.  Now the center $c$ acts as zero on $W'$. The
 $\Vir[G]=\Vir[b_1',G_0]$-module $V$ is generated by $W'$ and
 $d_{k b_1'+a_0}W'=0$ for any $k \in \N, a_0\in G_0$.
 So $V$ is the unique irreducible quotient of
 $M(b_1',G_0,W')$. If $W'=\bC
 v_0$ then $V=\bC v_0$. Since $V$ is nontrivial, we have $W'\simeq V'(\a,\b, G_0)$
 for some $(\a,\b)\in \bC^2$ and
 $V\simeq   M(b_1', G_0,
 V'(\a,\b, G_0))$.
\hfill $ $\qed
\enddemo

For any $\Z$-basis $B'=\{b_1',b_2',\cdots  ,b_n'\}$ of $G$, we
define the total order "$\su_{B'}$"  on G as follows:
$x_1b_1'+x_1b_2'+\cdots  +x_nb_n'\su_{B'}y_1b_1'+y_1b_2'+\cdots
+y_nb_n'$ if   $(x_1, x_2,\cdots  ,x_n)\su (y_1,y_2,\cdots
,y_n)$.

 \proclaim{Corollary 3.4} Suppose that $G\simeq \Z^2$.
 For any $(0,0)\neq (\dot{c},h) \in \bC^2$,
 and any $\Z$-basis $B'=\{b_1',b_2'\}$ of
 G, there exists $\lambda \in \supp (V(\dot{c},h,\su_{B'}))$ such
 that $\dim (V(\dot{c},h,\su_{B'}))_{\lambda}= \infty$.
\endproclaim
\demo{Proof} Suppose that for any $\lambda \in \supp
(V(\dot{c},h,\su_{B'}))$ we have  $\dim
(V(\dot{c},h,\su_{B'}))_{\lambda}< \infty$. It is easy to see that
$\supp (V(\dot{c},h,\su_{B'})) \subset (h-\N b_1'+\Z b_2')\bigcup
(h-\Z^+b_2')$, hence $V(\dot{c},h,\su_{B'})$ is a GHW module with
GHW $h$ w.r.t. $ B'$. Note that
 $$(h+ \N b_1'+\Z b_2')\bigcap \break
\supp (V(\dot{c},h,\su_{B'}))=\emptyset,\,\,\,\,\roman{and} \,\,\,
$$ $$(h+\Z b_2')\bigcap \supp (V(\dot{c},h,\su_{B'}))\neq
\emptyset.$$ Using the same argument as in the proof of Claim 2 of
Lemma 3.2,  we see that $W=\bigoplus_{\lambda \in \Z b_2'}
V(\dot{c},h, \su_{B'})_{h+\lambda}$ is a uniformly bounded
$\Vir[b_2']$ module. Since W contains the submodule
$W'=U(\Vir[b_2'])(v_h)$  which is a highest weight module with
highest weight $(\dot c,h)$, $W'$ (and $W$) is not uniformly
bounded. A contradiction. Hence $(\dot c,h)=(0,0)$. The corollary
follows.
 \hfill $ $\qed
\enddemo

 \par
 \proclaim{Lemma 3.5}  Suppose that $G\simeq \Z^2$. If there exist
  $(k,l) \neq 0,\,\,(i,j) \in\Z^2$,  $p,q \in \Z$
 such that$$\{ x \in \Z|\Lambda_0+ib_1+jb_2+x(kb_1+lb_2)
 \in \supp(V)\}\supset(-\infty, p] \bigcup [q, \infty),$$
 then   $V\simeq   M(b_1', \Z b_2', V'(\a,\b,\Z b'_2))$ for some $ \a,\b \in
\bC$, and a $\Z$-basis $B'=(b_1',b_2')$ of G.

\endproclaim
\demo{Proof }
  From Lemma 3.1 d6), we see $kl<0.$ We may assume  that $l>0$.
Let $$(i_0,j_0)=\cases(i+qk,j+pl), &\,\,\,if\,\,\,p<q-1,\\
(i,j), &\,\,\,if\,\,\,p\geq q-1.\endcases$$
  Denote $L:=\{i_0b_1+j_0b_2+x(kb_1+lb_2)|x\in\Z\}$.
If $p<q-1$, write $i_0b_1+j_0b_2+x(kb_1+lb_2)=
ib_1+jb_2+(x+q)(kb_1+lb_2)+(p-q)l_0b_2$ or
$ib_1+jb_2+(x+p)(kb_1+lb_2)+(q-p)k_0b_1$ according to $x\ge0$ or
$x<0$.
  From Lemma 3.1 d5) we see that all points in the set
   $\Lambda_0 +L$ are weights of $V$.

  Write $(k,l)=s(k_0,l_0)$ with $k_0,l_0$ relatively prime, $s\geq 1$.
  By replacing $(i,j)$ with $(i_0,j_0-(s-1)l_0)$, we may assume that
  $p=q$. Then similarly we have
$L_0$:$=ib_1+jb_2+\Z(k_0b_1+l_0b_2)\subset
  \supp (V)$. Using Lemma 3.1 d5) we see
that $$\{ \Lambda_0+x{b_1}+y {b_2} | l_0x-k_0y \leq
l_0i-k_0(j_0-(s-1)l_0), (x,y)\in \Z^2 \} \subset
  \supp (V),$$
  i.e., all points under the line  $\Lambda_0 -(s-1)l_0b_2+L_0$
are weights of $V$ (It might help if one draws a diagram on the
$Ob_1b_2$-plane).

  So we may assume
  that $k,l$ are relatively prime, $k<0$, $l>0$,
  and there exists an integer $m_0$ such that

  $$\{ \Lambda_0+x{b_1}+y {b_2} | lx-ky \leq m_0, (x,y)\in \Z^2 \} \subset
  \supp (V). \tag 3.7$$

  Fix $(k',l')\in \Z^2$ with
  $lk'-kl'=1$.  Denote $b_1'=kb_1+lb_2$ and  $b_2'=k'b_1+l'b_2$.
If $$ \{\Lambda_0-kb_1+b_2'+tb_1'|t \in \Z\}\bigcap \supp (V)=
 \emptyset,$$
 then the lemma follows from Lemma 3.3. Hence we may assume that
 $$ \{\Lambda_0-kb_1+b_2'+tb_1'|t \in \Z\}\bigcap \supp (V)\neq
 \emptyset. \tag 3.8$$
Choose $\Lambda_0-kb_1+b_2'-sb_1' \in \supp(V)$, and a nonzero
weight vector $v\in V_{\Lambda_0-kb_1+b_2'-sb_1'} $. Let
$b_1''=sb_1'-b_2'$, $b_2''=(s+1)b_1'-b_2'$. Since
$\Lambda_0-kb_1,\Lambda_0-kb_1+b_1'\in \Lambda_0+\Z^+b_1+\Z^+b_2$,
we obtain
$$d_{b_1''}v=0,\,\,\,\,
d_{b_2''}v=0.$$ Thus
   $$ d_{mb_1''+nb_2''}v=0 \,\,\forall\,\, m>0,
   \,\,n>0.$$
    Using this, one sees that $v$ is a GHW vector
    with respect to the $\Z$-basis $\{b_1''+b_2'',b_1''+2b_2''\}$ of G. Now
    by Lemma 3.1 b)  there exists some $x_0$ such that
    $$\lambda_0+b_2'+x((b_1''+b_2'')+(b_1''+2b_2''))\notin \supp (V),\,\,\,\forall\,\,\,
    x>x_0.\tag3.9$$
    But
    $$(b_1''+b_2'')+(b_1''+2b_2'')=2b_1''+3b_2''
    =(2s+3(s+1))b_1'-5b_2'$$ $$\hskip 3cm
    =((5s+3)k-5k')b_1+((5s+3)l-5l')b_2,$$
    $$l((5s+3)k-5k')-k((5s+3)l-5l')=-5(lk'-kl')=-5<0. $$ Hence
    for $x$ sufficiently large we have
  $$\lambda_0+b_2'+x((b_1''+b_2'')+(b_1''+2b_2''))\in \{
  \Lambda_0+x{b_1}+y {b_2} | lx-ky \leq m_0, (x,y)\in \Z^2
  \},
  $$ which is  a contradiction to (3.7) and
  (3.9), hence (3.8) cannot occur.  The lemma follows.
 \hfill $ $\qed
\enddemo

 \proclaim{Lemma 3.6}  Suppose that $G\simeq \Z^2$. If there exist  $(i,j), (k,l) \in \Z^2$
 and $x_1,x_2,x_3 \in \Z$ with
 $x_1<x_2<x_3$, such that
 $$ \Lambda_0+ib_1+jb_2+x_1(kb_1+lb_2)\notin \supp (V),$$
  $$\Lambda_0+ib_1+jb_2+x_2(kb_1+lb_2) \in \supp (V), \,\,\,and$$
  $$\Lambda_0+ib_1+jb_2+x_3(kb_1+lb_2) \notin \supp (V),$$
  then a). there exists  $x\in \Z$  with $x_1<x<x_3$ such that
$$\Lambda_0+ib_1+jb_2+x(kb_1+lb_2)=0,$$ and  further,
 b). such a module $V$ does not exist.
  \endproclaim

\demo{Proof} We may assume that $k,l$ are relatively prime, and by
Lemma 3.1 d6) we see $kl<0$. So we may assume that $k<0$ and
$l>0$. Replacing $x_2$ by the largest $x<x_3$ with
$\Lambda_0+ib_1+jb_2+x(kb_1+lb_2) \in \supp (V)$, and then
replacing $x_3$ by $x_2+1$ and $(i,j)$ by $(i,j)+x_2(k,l)$ we can
assume that
$$
x_1<x_2=0,\,\,\,x_3=1.\tag 3.10
$$
 Fix a nonzero weight vector $v  \in V_{\Lambda_0+ib_1+jb_2}$.
 Then (3.10) means
 $$
 d_{kb_1+lb_2}v =0=d_{x_1(kb_1+lb_2)}v ,
 $$
 which yields $d_{\pm(kb_1+lb_2)}v =0$.
 By Lemma 3.1 b) we can choose $p,q> 0$ such that
 $d_{pb_1+qb_2}v =0$. Since $kq-lp <0$, then
 $S=\{b_1'=kb_1+lb_2,b_2'=pb_1+qb_2\}$ is a $\Z$-linear independent
 subset of G. Note that
 $d_{mb_1'+nb_2'}$  for $n>0$ belong to the subalgebra
  generated by $d_{\pm b_1'}, d_{b'_2}$. Thus
 $$d_{mb_1'+nb_2'}v =0,\,\,\forall \,\,n>0, m\in \Z.$$
 Consider the
 $\Vir[b_1']$-module $W=U(\Vir[b_1'])v  $.
 By using the PBW basis of $U(\Vir[S])$ we have
 $$(\Lambda_0+ib_1+jb_2+\Z b_1'+ \N b_2')\bigcap
 \supp (U(\Vir[S])v )=\emptyset. \tag 3.11$$

  {\it Case 1}:  $W$ is not uniformly bounded.

   From Virasoro algebra theory we see that $W$ has a nontrivial
   irreducible sub-quotient $\Vir[b_1']$-module
    $W_1/W_2$ which is a  highest (or lowest) weight $\Vir[b_1']$-module. Using  (3.11) and PBW Theorem,
    we know that    $W'=U(\Vir[S])W_1/U(\Vir[S])W_2$ is a highest weight $\Vir[S]$-module
    w.r.t. the lexicographic order determined by $\{b'_1,b_2'\}$ with
    highest weight not equal to (0,0). Now by Corollary 3.4, $W'$
    has a weight space of infinite dimension. So does $S$. This case does not
    occur.

 {\it Case 2}:  $W$ is uniformly bounded.

 First we can easily see that  the center $c$ acts as zero on $V$.
 From the fact that \break $\supp(V'(\a,\b,\Z b'_1))
 =\a+\Z b'_1$ or $\Z b'_1\setminus\{0\}$ and the assumption (3.10),
 we know that $W \subset V_0$, the weight space with $0$ weight. We deduce a).

  It is clear that $W=V_0=\bC v_0$. Denote by $W''$ the $\Vir[S]$-module
  generated by $W$, which is a $\Vir[S]$-submodule of
   $V$. Now by (3.11),  $\Vir[S]$-module $W''$ is a quotient module of
   $M(0,0,\su_{B'})$, and $W''$ is nontrivial (from Lemma 3.1 c)), so $W''$ is
reducible.  If $d_{-b_2'+s_0b_1'}v_0=0$ for some $s_0$, from
    $$d_{-b_2'+sb_1'}v_0
    =(-b_2'+(2s_0-s)b_1')^{-1}[d_{-b_2'+s_0b_1'},d_{s b_1'-s_0 b_1'}]v_0=0,$$ and the fact
    that
    $\{d_{-b_2'+sb_1'}| s \in \Z \}$  generates  $\{d_{-tb_2'+sb_1'}|
    s \in \Z, t \in
    \N \}$, combining with (3.11) we deduce that $W''$ is a trivial $\Vir[S]$-submodule,
    a contradiction to Lemma 3.1 c). So we
    have  $d_{-b_2'+sb_1'}v_0 \neq 0$ for any $s \in \Z$. Thus
     $\{-b_2'+sb_1'| s \in \Z\} \subset \supp(V).$
     Now by lemma 3.5, we have $V\simeq
  M(b_2', \Z b_1', V'(\a,\b,\Z b'_1))$ for some $ \a,\b \in
\bC$, and a $\Z$-basis $B'=(b_1',b_2')$ of G. It is easy to see
that $  M(b_2', \Z b_1', V'(\a,\b,\Z b'_1))$ does not satisfy
Condition a). Thus such a module $V$ does not exist.

This completes the proof.
 \hfill $ $\qed
\enddemo

The idea of Claims 1 and 2 in the  proof of the next theorem comes
from the proof of [S2, Theorem 1.1] for $n=2$.

\proclaim{Theorem 3.7} Suppose that $B=({b_1},{b_2})$ is a
$\Z$-basis of the additive subgroup $G\subset \bC$. If $V$ is a
nontrivial
 irreducible weight module with finite
dimensional weight spaces over the higher rank Virasoro algebra
$\Vir[G]$, then $V\cong V'(\a,\b,G)$ or $V\cong   M(b_1', \Z b_2',
V'(\a,\b,\Z b_2' ))$ for some $ \a,\b \in \bC$, and a $\Z$-basis
 $B'=(b_1',b_2')$ of G.

\endproclaim
\demo{Proof } To the contrary, we suppose  that $V\ncong
V'(\a,\b,G)$ or $ M(b_1', \Z b_2', V'(\a,\b, \Z b'_2))$ for any $
\a,\b \in \bC$, and any $\Z$-basis of $B'=(b_1',b_2')$ of $G$.
From Theorem 2.5 we may assume that $V$ is a GHW module with GHW
$\Lambda_0$ w.r.t. the basis $B=\{b_1,b_2\}$ for $G$.  We need to
prove that  $V\cong   M(b_1', \Z b_2', V'(\a,\b,\Z b_2' ))$ for
proper parameters. We still assume that $B$ satisfies Lemma 3.1
d). By Lemmas 3.3, 3.5 and 3.6, for any $(i,j), 0\neq (k,l) \in
\Z^2$, there exists
 $p \in \Z $ such that
 $$\{x \in\Z|\Lambda_0+ib_1+jb_2+x(kb_1+lb_2) \in \supp(V)\}=(-\infty,p]
 \,\,or\,\,[p, \infty). \tag3.12$$

 Then for any  $i \in \N,$ there exist $x_i,y_i\in\Z^+$ such that
 $$(-\infty,y_i]=\max\{y \in \Z|\Lambda_0-ib_1+yb_2 \in \supp(V)\}, $$
 $$(-\infty,x_i]=\max\{x \in \Z|\Lambda_0+xb_1-ib_2 \in \supp(V)\}.$$
By lemma 3.1 d5) we know that $y_{i+1}\geq y_i \geq 0, x_{i+1}\geq
x_i\geq 0$. Let $j, t \in \N$, if $y_{jt} \geq t(y_j+1)$, then
$t>1$ and $\Lambda_0, \Lambda_0+t(-jb_1+(y_j+1)b_2)\in \supp(V)$,
and  by (3.12), $\Lambda_0+(-jb_1+(y_j+1)b_2) \in \supp(V)$,
contrary to the definition of $y_j$. So
$$ y_{tj} < t(y_j+1),\,\,
\forall \,\,t, j \in \N. \tag 3.13
$$
Since $\Lambda_0+b_2 \notin \supp(V)$ and
$\Lambda_0-jb_1+y_jb_2=\Lambda_0+b_2+(-jb_1+(y_j-1)b_2) \in
\supp(V)$, then (3.12) yields $\Lambda_0+b_2+t(-jb_1+(y_j-1)b_2)
\in \supp(V)$ for all $t>0$. Hence
$$y_{tj} \geq t(y_j-1)+1, \,\,\forall \,\, t, j
\in \N. \tag 3.14
$$
Using (3.13), (3.14) we obtain $$j(y_i-1)+1 \leq y_{ij} <i(y_j+1)
\,\,\forall \,\, i,j \in \N.$$ From $j(y_i-1)+1 <i(y_j+1)$ and the
 one with  $i,j$ interchanged,  we deduce
$$\frac{y_j}{j}-\frac{i+j-1}{ij}<\frac{y_i}{i}
<\frac{y_j}{j}+\frac{i+j-1}{ij},\,\,\forall \,\, i,j \in
\N.\tag3.15$$  This shows that the following limits exist:
$$\a=\lim_{i\rightarrow \infty} \frac{y_i}{i},\,\,
 \b=\lim_{i\rightarrow \infty} \frac{x_i}{i},\tag3.16$$
 where the second equation is obtained by symmetry.
 Note that (3.12) implies that there exists some $j_0\in \N$ such that $y_{j_0}>1$
 (otherwise $(\Lambda_0+2b_2+\Z b_1)\cap \supp(V)=\emptyset$).
 Hence by (3.14) we deduce
$$\frac{y_{tj_0}}{tj_0} \geq \frac{y_{j_0}-1}{j_0}+ \frac{1}{tj_0}>
\frac{y_{j_0}-1}{j_0}>0,$$ thus $\a>0$. Similarly $\b>0$.

\vskip 5pt
 {\bf Claim 1.} {\it $\a=\b^{-1}$ is an irrational number.}

\vskip 5pt
 {\it Proof of Claim 1.}
Suppose that $\a>\b^{-1}$. Choose $s,q \in \N$ with $s,q$
relatively prime and $\a>\frac{s}{q}>\b^{-1}$. Applying (3.12) to
$\Lambda_0+t(-qb_1+sb_2)$, by the definition of $\a$,  we have
$\Lambda_0+t(-qb_1+sb_2) \in \supp(V)$ for all sufficiently large
$t$. From (3.12), hence, for all sufficiently large $t$ we have
$\Lambda_0-t(-qb_1+sb_2) \notin \supp(V)$, which implies that
$\b=\lim_{t\to\infty}\frac{x_{st}}{st}\le \frac{q}{s}$, i.e ,
$\b^{-1}\ge \frac{s}{q}$, a contradiction. So we have $\a \leq
\b^{-1}$, and similarly we have $\a \geq \b^{-1}$. Thus
$\a={\b}^{-1}$.

 Assume $\a= \frac{q}{s}$ is a rational number, where $s,q \in \N$
 are
 relatively prime. By (3.12), there exists some $m_0 \in \Z$  such that
 $\Lambda_0-b_2+m_0(sb_1-qb_2)\notin \supp(V)$. Say $m_0>0$.
 Since $\Lambda_0 \in \supp(V)$, by
 (3.12) again, we deduce $$\Lambda_0+i(-m_0sb_1+(m_0q+1)b_2) \in \supp(V),\,\,\forall\,\,i
 \in [0, \infty].$$
However
$\a=q/s=\lim_{i\to\infty}\frac{y_{im_0s}}{im_0s}\ge\frac{_{m_0q+1}}{m_0s}>q/s$,
a contradiction. Hence $\a$ is an irrational number, and Claim 1
follows.
 \qed

\vskip 5pt
 We define a total order $>_{\a}$ on G as follows:
$$ib_1+jb_2>_{\a} kb_1+lb_2 \Longleftrightarrow i\a+j>k\a+l.$$
 Let $G^+\hskip - 2pt =\{ib_1+jb_2\in G |ib_1+jb_2>_{\a}0\}$. If
 $\lambda\in\supp(V)$ satisfies
 $(\lambda+G^+)\bigcap \supp(V)\break =\emptyset$, then $V$ is a nontrivial highest
 weight module w.r.t. ``$<_{\a}$". Since the order ``$<_{\a}$" is dense,
 from Theorem 2.2 we see that $V$ is a Verma module, which contradicts the
 fact that
 all weight spaces of $V$ are finite dimensional.
 So for any $\lambda\in\supp(V)$ we have $$(\lambda+G^+)\bigcap \supp(V)\neq \emptyset.\tag 3.17$$

\vskip 5pt
 {\bf Claim 2.} {\it If $\Lambda_0+g \in \supp(V)$ for some $g=ib_1+jb_2 \in
 G^+$ then
  $$\Lambda_0+kb_1+lb_2 \in
  \supp(V)\,\,\,\,\forall\,\,  kb_1+lb_2<_{\a}ib_1+jb_2. $$}

\vskip 5pt
 {\it Proof of Claim 2.}
  If there exists some
 $kb_1+lb_2<_{\a}ib_1+jb_2$ such that $\Lambda_0+kb_1+lb_2 \notin
  \supp(V),$ (3.12) implies
  $$\Lambda_0+ib_1+jb_2+t((k-i)b_1+(l-j)b_2) \notin \supp(V), \,\,\forall \,\, t \in
  \N.$$ If $k-i<0$ (then $l-j\ge0$), from $(k-i)b_1+(l-j)b_2<_{\a}0$ we see that $-(l-j)/(k-i)<
  \a$. On the other hand,
  $$\a=\lim_{t\to\infty}\frac{y_{t(i-k)-i}}{t(i-k)-i}
  \le \lim_{t\to\infty}\frac{j+t(l-j)}{t(i-k)-i}
  =\frac{{l-j}}{i-k},$$
   a contradiction.
  If $k-i>0$,
   from $(k-i)b_1+(l-j)b_2<_{\a}0$ we know that $(l-j)<0$,
  and  $-(k-i)/(l-j)<\a^{-1}$. Similarly,
  $\a^{-1}=\lim_{t\rightarrow \infty} \frac{x_t}{t}\le -(k-i)/(l-j)$,
  again a contradiction.
  If $ k-i=0$, by Lemma 3.1 d5) we have $(l-j)>0$, but by
  $(k-i)b_1+(l-j)b_2<_{\a}0$, we have $l-j<0$, which is also a contradiction.
  So we have Claim 2. \qed
\vskip 5pt Claim 2 implies that for any $\Lambda\in\supp(V)$, we
have
$$\Lambda-G^+\subset \supp(V).\tag3.18$$

\vskip 5pt
 {\bf Claim 3.} {\it $d_{-g} v_{\lambda} \neq 0$ for any
$g=ib_1+jb_2 \in G^+$ and any nonzero weight vector $
v_{\lambda}\in V_{\lambda}$.} \vskip 5pt
 {\it Proof of Claim 3.}
 Suppose that $d_{-g} v_{\lambda} = 0$ for some $g=ib_1+jb_2 \in G^+$ and
$0\neq v_{\lambda}\in V_{\lambda}$. By  (3.12) and (3.18), we see
that $d_{sg}v_{\lambda}=0$ for all sufficiently large  $s> 0$.
Hence $d_{g}v_{\lambda}=0$. By Lemma 3.1 b) we can choose
$g_1=pb_1+qb_2$ such that
 $d_{g_1}v_{\lambda}=0$ and $S=\{g,g_1\}$  is a $\Z$-linearly independent
 subset of G. Consider the $\Vir[g]$-module W=
 $U(\Vir[g]) v_{\lambda}$. Using the PBW basis of $U(\Vir[S])$ we have
 $$(\lambda+\Z g+ \N g_1)\bigcap
 \supp (U(\Vir[S])v_{\lambda}=\emptyset. $$
 By  (3.12) and (3.18) there exists some $s_0$ such that  $\lambda+sg
 \notin \supp(V)$ for all $s>s_0$. Hence any irreducible $\Vir[g]$-subquotient  of $W$ is
 a  highest weight module. If $W$ has a nontrivial irreducible $\Vir[g]$-subquotient,
 using the arguments, analogous to those
used in Case 1 in the proof of
 Lemma 3.6, we get a contradiction. So we deduce that $W= \bC v_\lambda$ with $\lambda=0$.
 With a similar discussion as  in Case
 2 in the proof of Lemma 3.6 we obtain that $\lambda+\Z g-g_1 \subset \supp
 (U(\Vir[S])v_{\lambda})$, which contradicts (3.12). Hence  Claim 3 follows. \qed

 \medskip

Fix  $\Lambda_0+ib_1+jb_2 \in \supp(V)$, where $ib_1+jb_2 \in
G^+$. We are going to show that $\dim
V_{\Lambda_0+ib_1+jb_2}=\infty$. For a given $n>0$, let
$\varepsilon=\frac{1}{n}(j+i\a)>0$. Since the order $``<_{\a}"$ is
dense, we can choose $p,q \in \Z$ with $0<q+p\a <\varepsilon$.
 Hence we obtain
 $0<_{\a}pb_1+qb_2$ and
 $npb_1+nqb_2<_{\a}(ib_1+jb_2)$. Then from Claim 2 we deduce that
 $$\Lambda_0+m(pb_1+qb_2) \in \supp(V)\,\,\,\forall\,\,\,m \leq 0.$$   By (3.12) we
 assume that $m_0$ is the maximal integer
 such that $\Lambda_0+m_0(pb_1+qb_2) \in \supp(V)$, so  $m_0\geq
 n$. Let $$M=\{g \in G^+ |0\ne \Lambda_0+m_0(pb_1+qb_2)+g \in \supp(V)\}.$$
 By (3.17) $M$ is an infinite set.
 Denote $\overline{g}=pb_1+qb_2$.

\vskip 5pt
 {\bf Claim 4.} {\it There exist  $g_0\in M$ such that for any $k:\,\,1\leq k \leq n$,
  the $k$ vectors $$d_{-\overline{g}}^{k-1} d_{-\overline{g}}v,\,\,\,d_{-\overline{g}}^{k-2} d_{-2\overline{g}}v ,\,\,\,\cdots  ,\,\,\,
d_{-\overline{g}}
d_{-(k-1)\overline{g}}v,\,\,\,d_{-k\overline{g}}v $$  are linearly
 independent, where $v\in
 V_{\Lambda_0+g_0+m_0\overline{g}}\setminus\{0\}$.}

\vskip 5pt
 {\it Proof of Claim 4.}
 We will prove the Claim by induction on $k$.

Suppose that $v\in V_{\Lambda_0+g+m_0\overline{g}}\setminus\{0\}$
for $g\in M$.

 If $k=1$, from $d_{\overline{g}}v_{\Lambda_0+g+m_0\overline{g}}=0$,  we
 deduce that
 $$d_{\overline{g}}d_{-\overline{g}}v
 =(-2\overline{g}(\Lambda_0+g+m_0\overline{g})+\frac{\overline{g}^3-\overline{g}}{12}c)v .$$
 Let $h_1(g):=-2\overline{g}(\Lambda_0+g+m_0\overline{g})+\frac{\overline{g}^3-\overline{g}}{12}c$.
Then the set $M_1=\{g\in M|h_1(g)\ne0\}$ is infinite and
$d_{-\overline{g}}v\ne0$ for any $g\in M_1$.

Suppose that $k>1$ and there exist an infinite set $M_{k-1}\subset
M$ such that \break $d_{-\overline{g}}^{k-2}
d_{-\overline{g}}v,\,\,\,d_{-\overline{g}}^{k-3}
d_{-2\overline{g}}v ,\,\,\,\cdots  ,\,\,\, d_{-\overline{g}}
d_{-(k-2)\overline{g}}v,\,\,\,d_{-(k-1)\overline{g}}v $  are
linearly
 independent for any $v\in
 V_{\Lambda_0+g+m_0\overline{g}}\setminus\{0\}$ and $g\in
 M_{k-1}$.

Now we consider $k$.
 If the vectors
$$d_{-\overline{g}}^k v ,\,\,
  d_{-\overline{g}}^{k-2} d_{-2\overline{g}}v ,
  \cdots  ,  d_{-\overline{g}}
d_{-(k-1)\overline{g}}v ,\,\,
 d_{-k\overline{g}}v $$ are linearly
 dependent for some $g \in M_{k-1}$.
 Then there exist  $a_1,\cdots  ,a_k\in\bC$, not all zero, such
 that
 $$w_{k}=a_1d_{-\overline{g}}^{k} v
 +a_2d_{-\overline{g}}^{k-2} d_{-2\overline{g}}v +\cdots+
 a_{k}d_{-(k)\overline{g}}v =0.$$
 Using $[d_{\overline{g}}, d_{-\overline{g}}^k]=
 -k\overline{g}(2d_0+(k-1)\overline{g}-\frac{\overline{g}^2-1}{12}c)d_{-\overline{g}}^{k-1}$,
 we deduce that

  $0=d_{\overline{g}}w_{k}$

  $=-a_1k\overline{g}(2(\Lambda_0+g+(m_0-k+1)\overline{g})
  +(k-1)\overline{g}-\frac{\overline{g}^2-1}{12}c)d_{-\overline{g}}^{k-1}v $

  ${ }\hskip .1cm +a_2(-k+2)\overline{g}(2(\Lambda_0+g+(m_0-k+1)\overline{g})+
  (k-3)\overline{g}-\frac{\overline{g}^2-1}{12}c)
  d_{-\overline{g}}^{k-3}d_{-2\overline{g}}v $

  ${ }\hskip .1cm -3a_2\overline{g}d_{-\overline{g}}^{k-1}v +\cdots  $

 ${ }\hskip .1cm +a_{k-1}(-1)\overline{g}(2(\Lambda_0+g+(m_0-k+1)\overline{g})-
 \frac{\overline{g}^2-1}{12}c)d_{-(k-1)\overline{g}}v $

   ${ }\hskip .1cm +(-k)a_{k-1}\overline{g}d_{-\overline{g}}d_{-(n-2)\overline{g}}
  v_{\Lambda_0+g+(m_0-k+1)\overline{g}}$

 ${ }\hskip .1cm +a_{k}(-k-1)\overline{g}d_{-(k-1)\overline{g}}v $.
  \newline
  This together with the inductive hypothesis  yields that
    $$a_i=a_1 f_{i}(g),\,\, \forall\,\,\,i=1,2,\cdots  ,k ,\tag 3.19$$
   where $f_{i}(X)$ is a polynomial of degree $i-1$ in $X$.  Using (3.19)  and
   the following computations

 ${ }\hskip .1cm 0=d_{k\overline{g}}w_{k}$

 ${ }\hskip .3cm =a_1(-k-1)\overline{g}(-k)\overline{g}\cdots  (-3)\overline{g}
 (-2\overline{g}(\Lambda_0+g+m_0\overline{g})+\frac{\overline{g}^3-\overline{g}}{12}c)v $

 ${ }\hskip .5cm +a_2(-k-1)\overline{g}(-k)\overline{g}\cdots  (-4)\overline{g}
 (-4\overline{g}(\Lambda_0+g+m_0\overline{g})+\frac{(2\overline{g})^3-2\overline{g}}{12}c)v $

 ${  }\hskip .5cm +\cdots  $

  ${ }\hskip .5cm +a_{k-1}(-k-1)\overline{g}(-2(k-1)\overline{g}
 (\Lambda_0+g+m_0\overline{g})+\frac{((k-1)\overline{g})^3-(k-1)\overline{g}}{12}c)v $

 ${ }\hskip .5cm +a_{k}(-2k\overline{g}(\Lambda_0+g+m_0\overline{g})
 +\frac{(k\overline{g})^3-k\overline{g}}{12}c)v $

 ${ }\hskip .3cm =a_1h_k(g)v ,$
 \newline
 where $h_{k}(X)$ is a
 polynomial of degree $k$ in $X$. Then $M_k=\{g\in M_{k-1}|h_k(g)\ne0\}$ is infinite and
the vectors in Claim 4 are linearly independent for $g_0\in M_k$.
 Hence  Claim 4 follows.
 \qed

\vskip 5pt

 From Claim 4  we know that $\dim V_{\Lambda_0+g_0+(m_0-n)\overline{g}} \geq n$ for some $g_0 \in M$
 and for all $n\in\N$.
 Noting that $g_0,\bar g\in G^+$, by Claim 3 we deduce that $\dim V_{\Lambda_{0}} \geq n$  for all $n\in\N$.
 Hence  $\dim V_{\Lambda_{0}}= \infty$.
  This proves that (3.12) cannot occur, and the theorem follows.
\hfill $ $\qed
\enddemo

\proclaim{Lemma 3.8} Suppose that $G=\Z b_1'\bigoplus G_0$. Then
$\supp(  M(b_1', G_0, V'(\a,\b,G_0 ))=\supp(V'(\a,\b,G_0))\bigcup
(\a+G_0-\N b_1')$.
\endproclaim
\demo{Proof}  It is clear that  $$\supp(  M(b_1', G_0,
V'(\a,\b,G_0))\subset \supp(V'(\a,\b,G_0))\bigcup (\a+G_0-\N
b_1'),$$ $$\supp(V'(\a,\b,G_0)) \subset \supp(  M(b_1', G_0,
V'(\a,\b,G_0 )).$$ Suppose that there exists  $d=\a+g_0-s b_1'
\notin \supp(  M(b_1', G_0, V'(\a,\b,G_0 ))$, where $s>0$ and $g_0
\in G_0$.

Choose $\alpha+g_1\in \supp(V'(\a,\b,G_0))\setminus\{0\}$. Let
$d'=g_1-g_0+sb'$. We see that $\alpha+g_1\in\supp(V)$.
  Fix  $v \in V'(\a,\b,G_0)_{\alpha+g_1}$.
 Let W be  the $\Vir[d', b_1']$-submodule generated by $v$. Then we have
 an irreducible
 sub-quotient module $W'$ of $W$ with $\alpha+g_1\in \supp(W')$,
 and $\alpha+g_1\pm d' \notin \supp(W')$. We get a
 contradiction to Lemma 3.6.
This completes the proof of  the lemma.\hfill $ $\qed
\enddemo

 From the lemma above we see that $\supp(  M(b_1', G_0, V'(\a,\b,G_0 ))$
  equals either $\a-\Z^+b'_1+G_0$ or $(-\Z^+ b'_1+G_0)\setminus\{0\}.$
Finally we can handle the general case.

 \proclaim{Theorem 3.9} If $V$ is a nontrivial
 irreducible weight module with finite
dimensional weight spaces over the higher rank Virasoro algebra
$\Vir[G]$ for $G\simeq \Z^n$ ($n\ge2$), then $V\simeq V'(\a,\b,G)$
or
 $V\simeq   M(b_1', G_0, V'(\a,\b,G_0 ))$ for some $ \a,\b \in
\bC$,  $b_1'\in G\setminus\{0\}$, and a subgroup $G_0$ of $G$ with
$G=\Z b_1' \bigoplus G_0$.
\endproclaim
\demo{Proof} From Theorem 2.5 we may assume that $V$ is a
nontrivial irreducible GHW module with GHW $\Lambda_0$ w.r.t.
$B=\{{b_1},{b_2},\cdots  ,b_n\}$  over $\Vir[G]$, where $B$ is a
$\Z$-basis of the additive subgroup  $G$ of $\bC$, then we need to
prove that $V\simeq  \break  M(b_1', G_0, V'(\a,\b,G_0 ))$. We
still assume that $B$ satisfies Lemma 3.1 d).

We shall prove this  by induction on $n$. For $n=2$ this is
Theorem 3.7. Now suppose that the theorem holds for any $n \leq
N-1$ where $ N\geq 3$. We shall prove $V\simeq   M(b_1', G_0,
V'(\a,\b,G_0 ))$ for $n=N$.

\vskip 5pt
 If there exist $g \in G$ and a corank $1$ subgroup $G_0$ of $G$
 such that $(\Lambda_0+g +G_0) \bigcap
\supp(V)
 \subset \{0\} $, then  the theorem follows from Lemma 3.3.
 (Indeed,  If $(\Lambda_0+g +G_0) \bigcap \supp(V)=\{0\}$,
suppose that $G_0= \Z a_1+\cdots  +\Z a_{N-1}$. We may assume that
$\Lambda_0+g=0$. Then $(a_1+ \Z 2 a_1+\cdots  +\Z 2
a_{N-1})\bigcap \supp(V)=\emptyset$. Using Lemma 3.3 we have the
theorem).
 So we may assume that for any $g\in G$ and any corank $1$ subgroup $G_0$,
  $$(\Lambda_0+ g +G_0) \bigcap
\supp(V)
 \nsubseteq \{0\}. \tag 3.20$$ Hence the $\Vir[G_0]$ module $V_{\Lambda_0+ g
 +G_0}=\bigoplus_{x
 \in G_0} V_{\Lambda_0+x+g}$ has a nontrivial irreducible sub-quotient.
 By Lemma 3.8, Theorem 2.5 and the inductive hypothesis,
 for any corank $1$ subgroup $G_0$ and any
 $g\in G$ there exist a subgroup $G_{0,1}$ of $G_0$,
 $\lambda_0' \in \Lambda_0+g+G_0$ and $g_{0,1} \in G_0\setminus\{0\}$
  with $G_0=\Z g_{0,1}\bigoplus G_{0,1}$ such that
$$\lambda_0'+G_{0,1}-\N g_{0,1} \subset \supp(V). \tag 3.21$$
Note that some other elements in $\lambda'_0+G_0$ can also be in
$\supp(V)$. Next we are going to show that under the assumption
(3.21) such a module $V$ does not exist.

\vskip 5pt
 {\bf Claim 1.} {\it There are no
$\lambda_0 \in \supp(V)$, $t_0 \in \Z$,  $   g_0, g_1 \in
G\setminus\{0\}$ or subgroups $G'_1\subset G_0'\subset G$ with
$G=\Z g_0\oplus G_0'$  and $G'_0=\Z g_1\oplus G_1'$ satisfying
 $$\lambda_0-\Z^+g_1+G_1', \,\,\,\,\lambda_0+t_0g_1+\Z^+g_1+G_1' \subset \supp(V) .$$ }

(If $t_0\le0$, then $\lambda_0+G_0'\subset \supp(V) .$)
 \vskip 5pt
 {\it Proof of Claim 1}. Suppose that there exist
$\lambda_0 \in \supp(V)$, $t_0 \in \Z$,  $   g_0, g_1 \in
G\setminus\{0\}$ and subgroups $G'_1\subset G_0'\subset G$ with
$G=\Z g_0\oplus G_0'$ and $G'_0=\Z g_1\oplus G_1'$ satisfying
 $$\lambda_0-\Z^+g_1+G_1', \,\,\,\,\lambda_0+t_0g_1+\Z^+g_1+G_1' \subset \supp(V) .$$
Choose $0 \neq (k_1,\cdots  ,k_N) \in \Z^N$, $k_i$ relatively
prime, such that  $G_0'=\break\{\sum_{i=1}^N x_ib_i| \sum_{i=1}^N
k_ix_i=0 \}$.

If there exist $i,j$ such that $k_ik_j\le0$, then there exists
$b'\in G'_0\setminus\{0\}$, $b'\ge 0$ with
$\{x\lambda_0-g_1+xb'\in\supp(V)\}=(-\infty, m_0]$, a
contradiction to the assumption (consider whether $b'\in G'_1$).
Then $k_ik_j
>0$ for all $i,j \in [1,N]$. Hence we may assume that $k_i> 0$ for
  all $i \in [1,N]$. Let $g_0=\sum _{i=1}^N s_i^{(N)}b_i$.
  Since $G'_0\oplus \Z g_0=G$ we have $\sum_{i=1}^N
  s_i^{(N)}k_i =\pm 1$.  By replacing $g_0$ with $-g_0$ if necessary,
  we may assume that $\sum_{i=1}^N
  s_i^{(N)}k_i =1$.
  Choose a basis of $G_1'$, say $\{b_1',b_2',\cdots  ,b_{N-2}'\}$. Take
  $b_{N-1}'=g_1, b_{N}'=g_0$,
   then $B'=\{b_1',b_2',\cdots  ,b_N'\}$ is a basis of $G$.

   \medskip
   {\bf Subclaim } { \it For any $N_0>0$ there exists $m_0 \in \N$ such
   that $[m_0, \infty) \subset \break \{\sum_{i=1}^N k_ix_i| x_i \geq
   N_0,i=1,2,\cdots  ,N\}$.}
   \medskip
   {\it Proof of Subclaim.} Note that $\sum_{i=1}^N
  s_i^{(N)}k_i =1$.
   Choose $n_0 \in \N$ such that
    $n_0+s_i^{(N)}\geq 0$ for all $i$. Note that $k_{1} > 0$. Take
   $m_0= \sum_{i=1}^N k_i(N_0+k_{1}n_0)$. Noting  that
    $$m_0+tk_{1}=(\sum_{i=1}^N
    k_i(N_0+k_{1}n_0))+k_{1}t,\,\,\,\forall\,\,\,
   t>0,\,\,\,\,\roman{and}
\,\,\, $$ $$m_0+tk_{1}+i=(\sum_{i=1}^N k_i(N_0+k_{1}n_0+i
s_i^{(N)}))+tk_{1},\,\,\,\,\roman{for} \,\,\,0 \leq i <k_{1},$$
    we have proved the subclaim. \qed
\medskip
Denote $b'_{N-1}=g_1=\sum _{i=1}^N s_i^{(N-1)}b_i$.
 Choose $N_0\in\N$ such that $N_0+t_0s_i^{(N-1)}>0$ for all $i$,
 then choose $m_0$ for this $N_0$ as in the subclaim above.
   By the subclaim for any $m\geq m_0$,
   there exists $(x_1,x_2,\cdots  ,x_n) \geq (N_0,N_0,\cdots  ,N_0)$, such that
   $m=\sum_{i=1}^N k_ix_i$. Then using the choice of $(k_1,\cdots  ,k_N)$
   one can easily verify that
   $mb_N'- \sum_{i=1}^N x_ib_i \in G_0'$.
   Using this we can write  $\lambda\in \lambda_0-mb_N'+G'_0$ as
    $$\lambda=\lambda_0+h_0- \sum_{i=1}^N
    x_ib_i,\,\,\,\,
        \lambda=\lambda_0+h_0+t_0g_1- ((\sum_{i=1}^N
        x_ib_i)+t_0g_1),$$
        where $h_0 \in G'_0.$ Noting that
        $\sum_{i=1}^N x_ib_i , ((\sum_{i=1}^N x_ib_i)+t_0g_1) \in \sum_{i=1}^N \Z^+b_i$,
   and the fact that
   $\lambda_0+h_0 \in \supp(V)$ or $\lambda_0+h_0+t_0g_1\in \supp(V)$,
   using Lemma 3.1 d5) we deduce
   $$\lambda_0-m_0b_N'+ G_0'-\Z^+ b_N' \subset \supp(V). \tag 3.22$$

 Fix some $\lambda_0'\in \Lambda_0+(\sum_{i=1}^N\N b'_i)$  such that
 $$\lambda_0' ,\,\,\,\,
 \lambda_0' \pm b_i' ,\,\,\,\,\lambda_0' \pm b_i'- b_N' \in \Lambda_0+
 \sum_{i=1}^{N}\Z^+b_i \,\,\,\,\roman{for \,\, all}
\,\,\, i \in [1,N].\tag 3.23$$

 Applying (3.21) to $\lambda_0'$ and $G'_0$ (replace $G_0$ by $G'_0$), since $N>2$ we have  $i_0 \in
 [1,N-1]$ and $s_0 \in \Z$ such that
 $$ \lambda_0'+s_0b_{i_0}' \in \supp(V).\tag3.24$$

 Denote $b_i''=-s_0b_{i_0}'-b_i'$ for all $i \in
 [1,N] \backslash {i_0}$ and $b_{i_0}''=-(s_0+1)b_{i_0}'-b_N'$.
 Fix a nonzero $v_{\lambda_0'+s_0b_{i_0}'} \in V_{\lambda_0'+s_0b_{i_0}'}$.
 It is easy to see that $\{b_1'',\cdots  ,b_N''\}$ forms a $\Z$-basis of G. By
 (3.23) we have
 $$d_{b_i''} v_{\lambda_0'+s_0b_{i_0}'}=0 \,\,\,\,\roman{ for\,\,
 all \,}\,\,\,i\in [1,N].$$ So we
 have
 $$d_b v_{\lambda_0'+s_0b_{i_0}'}=0 \,\,\,\,\roman{ for\,\, all
 \,}\,\,\,b \in (\Z^+ b_1''+\Z^+ b_2''+\cdots  +\Z^+ b_N'')
 \setminus(\bigcup_{i=1}^N\Z^+ b''_i).\tag3.25$$
 Hence $v_{\lambda_0'+s_0b_{i_0}'}$ is a highest weight vector w.r.t.
$ B'''=\{ 2b_1''+b_2'',b_1''+b_2'', b_1''+b_3'',\cdots
,b_1''+b_N''\}$ which satisfies Lemma 3.1 d).
 Now by lemma 3.1 b) there exists $x_0 $ such that for any
 $x>x_0$,
$$
\lambda_0'+s_0b_{i_0}'+
x((2b_1''+b_2'')+(b_1''+b_2'')+(b_1''+b_3'')+\cdots
+(b_1''+b_N'')) \notin \supp(V).\tag3.26
$$
Write $(2b_1''+b_2'')+(b_1''+b_2'')+(b_1''+b_3'')+\cdots
+(b_1''+b_N'')= h_0-l'b_N'$,
$\lambda_0'=\lambda_0-m_0b_N'+g_0+lb_N'$ where $g_0,h_0\in G'_0$,
$l,l'\in\Z$, $l'>0$ (since $G'_0=\sum_{i=1}^{N-1}\Z b'_i$). Then
for sufficiently large $x$,
$$
\lambda_0'+s_0b_{i_0}'+ x((2b_1''+b_2'')+(b_1''+b_2'')+
(b_1''+b_3'')+\cdots  +(b_1''+b_N''))
$$
$$
\in\lambda_0-m_0b_N'+(l-l'x)b_N' +G'_0\subset \lambda_0-m_0b_N'+
G_0'-\N b_N',  \tag 3.27
$$
which contradicts (3.22). Thus Claim 1
follows. \qed
 \vskip 5pt

Denote  $\bar G_{t}=tb_1+\Z b_2+\Z b_3+\cdots  +\Z b_N $ for $t
\in \Z$.

\vskip 5pt

{\bf Claim 2}. {\it If for $\lambda_0 \in \Lambda_0+G $, $
g_1,g_1' \in \bar G_0 \setminus\{0\}$, and  subgroups of $\bar
G_0$: $G_1,G_1'$
  with $\bar G_0=\Z g_1\bigoplus G_1$, $\bar G_0=\Z g_1'\bigoplus
  G_1'$, we have
$$
\lambda_0-\N g_1 +G_1,\,\, \lambda_0-\N g_1' +G_1'
 \subset \supp(V),
$$
then $G_1=G_1'$.}

\vskip 5pt
{\it Proof of Claim 2}. Suppose that $G_1\ne G_1'$.
  Fix $0 \neq f_1=\sum_{i=1}^n u_ib_i \in \bar G_0$ (then $u_1=0$) satisfying
 $$f_1 \in -\N g_1+G_1\,\,\,\,\roman{and}
\,\,\, \sum_{i=1}^n u_ix_i=0
$$
 for all $(x_1,\cdots  ,x_N)\in\Z^n$ with
 $\sum_{i=1}^n x_ib_i \in G_1,$ and fix  $ 0 \neq f_1'
 =\sum_{i=1}^n u_i'b_i \in \bar G_0$ satisfying
 $$f'_1 \in -\N g_1'+G'_1\,\,\,\,\roman{and}
\,\,\, \sum_{i=1}^n u_i'x_i=0
$$
 for all $(x_1,\cdots  ,x_N)\in\Z^n$ with
$ \sum_{i=1}^n x_ib_i \in G_1'. $ (We simply write $f\perp G_1,
f'\perp G'_1$).
  Since $G_1\ne G_1'$ we see that $\Z f'_1\cap\Z f_1=\{0\}$. Hence we can
  choose a base $B'=\{b_1',b_2',\cdots  ,b_{N-1}'\}$ of
 $\bar G_0$ as follows:
 Fix $b_1'=\sum_{i=1}^N s_i^{(1)}b_i \in \bar G_0$ such that
  $s_1^{(1)},s_2^{(1)},\cdots,s_N^{(1)}$ are relatively prime,
 $$\sum_{i=1}^N u_is_i^{(1)}>0,\,\,\,\,\roman{and}
\,\,\, \sum_{i=1}^N
 u_i's_i^{(1)}<0,\tag3.28
 $$ and extend it to a $\Z$ basis
 $B'=\{b_1',b_2',\cdots  ,b_{N-1}'\}$
 of $\bar G_0$. By replacing $b_j'
 (j>1)$
 with $b_j'+mb_1', m\gg 0$ if necessary, we may assume that $b_j'=\sum_{i=1}^N
 s_i^{(j)}b_i$ satisfies
 $$\sum_{i=1}^N s_i^{(j)}u_i>0 \,\,\,\,\roman{and}
\,\,\, \sum_{i=1}^N
 s_i^{(j)}u_i'<0 \,\,\,\,\roman{for\,\, all}
\,\,\,j \in [1,N-1].\tag3.29$$ Since $f\perp G_1, f'\perp G'_1$,
we see that
$$b'_i\in-\N g_1+G_1\,\,\,\,\roman{and}
\,\,\, b'_i\in\N g'_1+G'_1, \,\,\,\,\roman{for\,\, all} \,\,\,
i\in[1,N-1].\tag 3.30$$
 Take $b_N'=b_1.$ Hence $B'=\{b_1',b_2',\cdots  ,b_N'\}$ is a basis of
 G. Fix
 $$\lambda_0'= \lambda_0+t_0b_N'+\bar g_0,\,\,\,
  \tag 3.31$$   where  $t_0>0, \bar g_0 \in \bar G_0$ are such that
 $$\lambda_0',\,\,\,\,\lambda_0' \pm b_i' ,\,\,\,\,\lambda_0'
 \pm b_i' - b_N' \in\Lambda_0 +
 \sum_{i=1}^{N}\Z^+b_i \,\,\,\,\roman{for\,\, all} \,\,\,i \in [1,N].$$
So
$$\lambda_0',\,\,\,\,\lambda_0' \pm b_i' ,\,\,\,\,\lambda_0' \pm b_i' - b_N' \notin
\supp(V) \,\,\,\,\roman{for\,\, all} \,\,\,i \in [1,N].\tag 3.32$$

 Now applying (3.21) to $\lambda_0'$ and $\bar G_0$, since $N>2$ we see that there exist some
$i_0\in[1,N-1]$ and $s_0 \in \N$ such that
$$ \lambda_0'+sb_{i_0}' \in \supp(V)\,\,\,\,\roman{for\,\, all}
\,\,\,s\geq s_0,\tag 3.33$$ $$
\roman{or}\,\,\,\lambda_0'+sb_{i_0}' \in \supp(V)
\,\,\,\,\roman{for\,\, all} \,\,\, s\leq -s_0. $$
 We may assume that (3.33) holds (if $\lambda_0'+sb_{i_0}' \in \supp(V)$ for
 all $s\leq -s_0$, then the remaining arguments are
 exactly the same, using $G_1$).
 Denote $b_i''=-s_0b_{i_0}'-b_i'$ for all $i \in
 [1,N]\backslash {i_0}$ and $b_{i_0}''=-(s_0+1)b_{i_0}'-b_N'$.

  From (3.30) we see that
$$ b_1''+\cdots  +b_N''\in-\sum_{i=1}^N b'_i-s_0nb'_{i_0}-\N
b_1\subset  -\N g_1'+G_1'-\N b_1.\tag 3.34$$

 Fix a nonzero $v_{\lambda_0'+s_0b_{i_0}'} \in V_{\lambda_0'+s_0b_{i_0}'}$.
 It is easy to see that $\{b_1'',\cdots  ,b_N''\}$ forms a $\Z$-basis of G and
 $$d_{b_i''} v_{\lambda_0'+s_0b_{i_0}'}=0\,\,\,\,\roman{for\,\, all}
 \,\,\, i\in [1,N].\tag3.35
 $$
 So we have
 $$
 d_b v_{\lambda_0'+s_0b_{i_0}'}=0 \,\,\,\,\roman{for\,\, all} \,\,\,
  b \in (\Z^+ b_1''+\Z^+ b_2''+\cdots  +\Z^+ b_N'')
 \setminus(\bigcup_{i=1}^N\Z^+ b''_i).\tag3.36
 $$
 Hence $v_{\lambda_0'+s_0b_{i_0}'}$ is a highest weight vector w.r.t.
$ B'''=\{ 2b_1''+b_2'',b_1''+b_2'', b_1''+b_3'',\cdots
,b_1''+b_N''\}$.
 Now by Lemma 3.1 b) there exists some $x_0 $ such that for any
 $x>x_0$ we have
$$\lambda_0'+s_0b_{i_0}'+ x(b_1''+b_2''+\cdots  +b_N'') \notin
\supp(V).\tag3.37$$

From (3.31) we can write $\lambda_0'=\lambda_0+t_0b_N'+lg'_1+h$
where $h\in G'_1$, $l\in\Z$. Using (3.34), for sufficiently large
$x$ we have
$$\lambda_0'+s_0b_{i_0}'+ x(b_1''+b_2''+\cdots  +b_N'') \in
 \lambda_0-\N g_1'
+G_1'-{\N}b_1\subset\supp(V)\tag3.38$$ since $\lambda_0-\N g_1'
+G_1'\in\supp(V)$.
 This is a
contradiction to (3.37). Hence $G_1=G_1'$ and Claim 2 follows.
\qed

\vskip 5pt
 Denote $ V_{\Lambda_0+\bar G_{t}}=\bigoplus_{ b \in \bar G_{t}}
V_{\Lambda_0+b}$ for $t\in\Z$. It is easy to see that
$V_{\Lambda_0+\bar G_{t}}$ is a $\Vir[b_2,\cdots  ,b_N]$-module.
For any $0 \neq \lambda \in \supp(V_{\Lambda_0+\bar G_{t}})$ (we
refer to (3.21) for the existence), $\lambda$ is a weight of a
nontrivial irreducible $\Vir[b_2,\cdots  ,b_N]$-subquotient  of
$V_{\Lambda_0+\bar G_{t}}$. From the inductive hypothesis and
Claim 1, we know that such a nontrivial irreducible
$\Vir[b_2,\cdots ,b_N]$-module is isomorphic to $ M(g_t, G_t,
V'(\a_t,\b_t,G_t))$ for suitable $ \a_t,\b_t \in \bC$, and $g_t,
G_t$ with $\bar G_0=\Z g_t \bigoplus G_t$. Thus from Lemma 3.8, if
$0\ne \lambda\in\supp(V)\cap (\Lambda_0+g_t+\bar G_0)$, then there
exists a corank $1$ subgroup $G_{\lambda}$ of $\bar G_0$ such that
$$\lambda+G_{\lambda} \subset \supp(V)\bigcup\{0\}.\tag3.39
$$
Combining this with Claims 1 and 2, we deduce that for any
$t\in\Z$ there exist a corank 1 subgroup $G_t$ in $\bar G_0$,
$\a_t\in \Lambda_0+\bar G_{t}$ and $g_t\in \bar G_0$ such that
$\bar G_0=\Z g_t \bigoplus G_t$ and
$$\supp(V_{\Lambda_0+\bar
G_{t}})\setminus\{0\}=( \a_t-\Z^+ g_t+G_t)\setminus\{0\}.\tag3.40
$$
In particular,
$$
\a_t-\N g_t+G_t\subset\supp(V_{\Lambda_0+\bar G_{t}}).\tag3.41
$$
Lemma 3.1 d5) and Lemma 3.8 ensure that
$$
\a_{t+1}-b_1-\N g_{t+1}+G_{t+1},\,\,\,\, \a_t-\N
g_t+G_t,\,\,\,\subset\supp(V_{\Lambda_0+\bar G_{t}}).
$$  It follows from Claim 2 that $G_t=
 G_{t+1}$ for all $t\in\Z$.
Thus there exist a corank 1 subgroup $G_0$ in $\bar G_0$, $\a_t\in
\Lambda_0+\bar G_{t}$ and $g_0\in \bar G_0$ with $\bar G_0=\Z^+
g_0 \bigoplus G_0$ such that either
$$\supp(V_{\Lambda_0+\bar
G_{t}})\setminus\{0\}=( \a_t+\Z^+
g_0+G_0)\setminus\{0\},\,\,\,\,\roman{or} \,\,\,
$$
$$\supp(V_{\Lambda_0+\bar
G_{t}})\setminus\{0\}=( \a_t-\Z^+ g_0+G_0)\setminus\{0\}.\tag3.42
$$

 If there exists $t \in \Z$ such that
$$\supp(V_{\Lambda_0+\bar
G_{t}})\setminus\{0\}=( \a_t-\Z^+ g_0+G_0)\setminus\{0\},
$$
$$\supp(V_{\Lambda_0+\bar
G_{t+1}})\setminus\{0\}=( \a_{t+1}+\Z^+ g_0+G_0)\setminus\{0\}.
$$
Similarly we have $\lambda_1,\lambda_2\in \Lambda_0+ \bar G_t$
such that
$$\lambda_1-\N g_0 +G_0,\,\, \lambda_2+\N g_0 +G_0
 \subset \supp(V),$$  which contradicts
Claim 1.
 So we may assume that
$$\supp(V_{\Lambda_0+\bar
G_{t}})\setminus\{0\}=( \a_t-\Z^+ g_0+G_0)\setminus\{0\}
,\,\,\,\forall \,\,\,t\in\Z.\tag 3.43
$$
Hence we may assume that $\alpha_t\in \Lambda_0+\Z g_0+\Z b_1$.
 Then for any
$\lambda\in\supp(V)$, we have
$$\lambda+G_0\subset\supp(V)\bigcup\{0\}.\tag 3.44$$

Consider the $\Vir[g_0,b_1]$-module $V_{\Lambda_0+\Z g_0+\Z b_1}$
where $g_0\in G_0\setminus\{0\}$ as before. From (3.43),
$V_{\Lambda_0+\Z g_0 + \Z b_1}$ has a nontrivial irreducible
$\Vir[g_0,b_1]$-subquotient (we refer to the last paragraph in the
proof of Lemma 3.2). Hence there exist some $\lambda_0'\in
\Lambda_0+\Z g_0 + \Z b_1$ and a basis $b_0',g_0' $ of $ \Z g_0+\Z
b_1$  such that
$$\lambda_0'+\Z g_0' \subset \supp(V_{\Lambda_0+\Z g_0 + \Z
b_1}).$$ From (3.43) with $t=0$ we know that $g_0' \notin \Z g_0.$
Hence by (3.44) $\lambda_0'+\Z g_0'+ G_{0} \subset \supp
(V)\bigcup \{ 0 \}$, and $\Z b'_0+(\Z g_0'+ G_{0})=G$, which
contradicts Claim 1. This completes the proof of the
theorem.\hfill $ $\qed
\enddemo

\vskip.3cm \Refs\nofrills{\bf REFERENCES}
\bigskip
\parindent=0.45in

\leftitem{[BB]} S. Berman and Y. Billig, Irreducible
representations for toroidal Lie algebras, {\it J. Algebra},
221(1999), no.1, 188-231.

\leftitem{[BZ]}  Y. Billig and K. Zhao,  Weight modules over
exp-polynomial Lie algebras, {\it J. Pure Appl. Algebra}, Vol.191,
23-42(2004).

\leftitem{[DZ]}  D.Z. Djokovic and K. Zhao, Derivations,
isomorphisms, and second cohomology of {generalized} Witt
algebras, {\it Trans.  Amer. Math. Soc.} (2) {\bf 350}(1998),
643-664.

\leftitem{[HWZ]} J. Hu, X. Wang, and K. Zhao, Verma modules over
generalized Virasoro algebras $\Vir[G]$, {\it J. Pure Appl.
Algebra}, 177(2003), no.1, 61-69.


\leftitem{[K]} I. Kaplansky, Seminar on simple Lie algebras, {\it
Bull. Amer. Math. Soc.}, 60(1954), 470-471.

\leftitem{[KR]} V. G. Kac and K. A. Raina, ``Bombay lectures on
highest weight representations of infinite dimensional Lie
algebras,'' World Sci., Singapore, 1987.

\leftitem{[M]} O. Mathieu, Classification of Harish-Chandra
modules over the Virasoro algebra, {\it Invent. Math.} {\bf
107}(1992), 225-234.

\leftitem{[Ma1]} V. Mazorchuk, Classification of simple
Harish-Chandra modules over ${\bQ}$-Virasoro algebra, {\it Math.
Nachr.} {\bf 209}(2000), 171-177.

\leftitem{[Ma2]} V. Mazorchuk,  Verma modules over generalized
Witt algebras, {\it Compositio Math.} (1) {\bf 115}(1999), 21-35.

\leftitem{[O]} J. M. Osborn, New simple infinite-dimensional Lie
algebras of characteristic 0, {\it J. Algebra} {\bf 185}(1996),
820-835.

\leftitem{[P]} D. S. Passman, Simple Lie algebras of Witt type,
{\it J. Algebra} {\bf 206}(1998), 682-692.

\leftitem{[PZ]} J. Patera and H. Zassenhaus, The higher rank
Virasoro algebras, {\it Comm. Math. Phys.} {\bf 136}(1991), 1-14.

\leftitem{[R]} R. Ree, On generalized Witt algebras, {\it Trans.
Amer. Math. Soc.} {\bf 83}(1956), 510-546.

\leftitem{[St]} H. Strade, Representations of the Witt algebra,
{\it J. Algebra} {\bf 49}(2)(1977), 595-605.

 \leftitem{[S1]} Y.  Su, Simple modules over the high rank Virasoro algebras,
 {\it  Comm. Algebra}, 29(2001), no.5, 2067-2080.

 \leftitem{[S2]} Y.  Su,
Classification of Harish-Chandra modules over the higher rank
Virasoro  algebras, {\it Comm. Math. Phys.},   {\bf 240} (2003),
539-551.

\leftitem{[SZ]} Y. Su and K. Zhao, Generalized Virasoro and
super-Virasoro algebras and modules of intermediate series, {\it
J. Algebra}, 252(2002), no.1, 1-19.

 \leftitem{[W]} R. L. Wilson, Classification of generalized
Witt algebras over algebraically fields, {\it Trans. Amer. Math.
Soc.} {\bf 153}(1971), 191-210.

\leftitem{[X]} X. Xu, New generalized simple Lie algebras of
Cartan type over a field with characteristic $0$, {\it J.
Algebra}, 224(2000), no.1, 23-58.

 \leftitem{[Z]} H. Zassenhaus, Ueber Lie'sche Ringe mit
Primzahlcharakteristik, {\it Hamb. Abh.} {\bf 13}(1939), 1-100.

\leftitem{[Zh]} K. Zhao, Weight modules over generalized Witt
algebras with $1$-dimensional weight spaces, {\it Forum Math.},
Vol.16, No.5, 725-748(2004).

\endRefs
\vfill
\enddocument
\end